\documentclass[leqno, 11pt]{amsart}
\usepackage[utf8]{inputenc}

\usepackage{iftex}
\ifpdf
\else
  \PassOptionsToPackage{dvipdfmx}{hyperref}
  \PassOptionsToPackage{dvipdfmx}{graphicx}
  \PassOptionsToPackage{dvipdfmx}{xcolor} 
  \PassOptionsToPackage{dvipdfmx}{xy}  
\fi

\usepackage{amsfonts,amsmath,amsthm,amssymb,amscd}
\usepackage{mathrsfs}

\usepackage{appendix}

\usepackage[all]{xy}
\usepackage[T1]{fontenc}
\usepackage{lmodern}

\usepackage{tikz}
\usetikzlibrary{intersections, calc, arrows, cd}
\usetikzlibrary{shapes.geometric}

\usepackage{colonequals}
\usepackage{fullpage}

\definecolor{cadmiumgreen}{rgb}{0.0, 0.42, 0.24}

\usepackage[
colorlinks, 
linkcolor=blue,
urlcolor=blue,
citecolor=cadmiumgreen,
pdfauthor={Omid Amini, Shu Kawaguchi}, 
pdftitle={},
pdfstartview ={FitV},
]{hyperref}

\newtheorem{Theorem}{Theorem}[section]
\newtheorem{Proposition}[Theorem]{Proposition}
\newtheorem{Lemma}[Theorem]{Lemma}
\newtheorem{Corollary}[Theorem]{Corollary}
\newtheorem{Claim}{Claim}[Theorem]

\newtheorem{Question}[Theorem]{Question}

\theoremstyle{definition}
\newtheorem{defii}[Theorem]{Definition}
\newenvironment{Definition}
  {\pushQED{\qed}\defii}
  {\popQED\enddefii}
\newtheorem{remm}[Theorem]{Remark}
\newenvironment{Remark}
  {\pushQED{\qed}\remm}
  {\popQED\endremm}
\newtheorem{exx}[Theorem]{Example}
\newenvironment{Example}
  {\pushQED{\qed}\exx}
  {\popQED\endexx}

\numberwithin{equation}{section}

\newcommand{\ZZ}{{\mathbb{Z}}}

\newcommand{\RR}{{\mathbb{R}}}
\newcommand{\CC}{{\mathbb{C}}}

\newcommand{\OO}{{\mathcal{O}}}

\newcommand{\Jac}{\operatorname{\operatorname{Jac}}}

\newcommand{\Div}{\operatorname{Div}}

\newcommand{\Pic}{\operatorname{Pic}}

\newcommand{\rank}{\mathbf r}

\newcommand{\rest}[2]{\left.{#1}\right\vert_{{#2}}}  
\renewcommand{\setminus}{\smallsetminus}

\newcommand{\wt}{\operatorname{wt}}

\newcommand{\Supp}{\operatorname{Supp}}

\newcommand{\val}{\operatorname{val}}
\newcommand{\Rat}{\operatorname{Rat}}

\newcommand{\WL}{\operatorname{WL}}
\newcommand{\WP}{\operatorname{WP}}
\newcommand{\slope}{\operatorname{sl}}
\newcommand{\zero}{\operatorname{div}}

\newcommand{\outdeg}{\operatorname{outdeg}}

\newcommand{\Proof}{{\sl Proof.}\quad}
\newcommand{\QED}{{\unskip\nobreak\hfil\penalty50\quad\null\nobreak\hfil
{$\Box$}\parfillskip0pt\finalhyphendemerits0\par\medskip}}

  {\end{list}}

  {\end{list}}

\begin{document}

\title{Weierstrass gap sequences and their weights on tropical curves}
\author{Omid Amini} 
\address{CNRS - Laboratoire de Math\'ematiques d'Orsay, 
Universit\'e Paris-Saclay, 91400 Orsay, France}
\email{omid.amini@universite-paris-saclay.fr}
\author{Shu Kawaguchi}
\address{Department of Mathematics, 
Kyoto University, Kyoto, 606-8502, Japan}
\email{kawaguch@math.kyoto-u.ac.jp}
\date{March 19th, 2026}


\maketitle

\begin{abstract}
Given a divisor on a tropical curve, we associate to each point of the curve a Weierstrass gap sequence. We investigate structural properties of these gap sequences and explore their relationship with the Weierstrass gap sequences of line bundles on algebraic curves via the tropicalization process.
\end{abstract}



%
\setcounter{Theorem}{0}

%
\setcounter{Theorem}{0}

\section{Introduction}
\label{sec:intro}
Let $\Gamma$ be a tropical curve of genus $g \geq 2$, and let $\rank=\rank_\Gamma$ denote the Baker--Norine rank function on the set of divisors on $\Gamma$ (see Section~\ref{sec:preliminaries} for an overview of divisor theory on tropical curves). Let $D$ be a divisor on $\Gamma$ of degree $d$ and nonnegative rank $\rank(D) = r$.  

To each point $p \in \Gamma$, we associate a strictly increasing sequence of positive integers 
\[
\underline{n} \colonequals (n_1, \ldots, n_{r+1}),
\]
called the \emph{Weierstrass gap sequence at $p$ with respect to $D$}, or simply the \emph{$D$-gap sequence of $p$}. An integer $m$ belongs to the $D$-gap sequence if 
\[
\rank(D-m(p)) <  \rank(D-(m-1)(p)).
\]
Since subtracting a point from a divisor decreases the rank by at most one, it follows that there are exactly $r+1$ gaps at $p$, and the largest gap is at most $d+1$ (see Section~\ref{sec:Sg:etc}). When it is necessary to emphasize $D$ and $p$, we write $\underline{n}_D(p)$ instead of $\underline{n}$. 

The weight of the $D$-gap sequence $\underline{n}$ is defined by 
\[
\wt(\underline{n}) \colonequals \sum_{i=1}^{r+1} (n_i - i) \in \ZZ_{\geq 0}.
\] 
A point $p \in \Gamma$ is called a {\em $D$-Weierstrass point} if $\rank(D- (r+1)(p)) \geq 0$. Since the sequence $\underline{n}$ is strictly increasing, this condition is equivalent to requiring that $\underline{n}_{D}(p) \neq (1, 2, \ldots, r+1)$. Denote by $\WL(D)$ the locus of 
$D$-Weierstrass points. This is a closed subset 
of $\Gamma$. 
For $p \in\WL(D)$, the $D$-Weierstrass weight $\wt_D(p)$ of $p$ is the weight of its $D$-gap sequence, that is,
$\wt_D(p) = \wt(\underline{n}_D(p))$. 

In the case where $D$ is the canonical divisor $K=K_\Gamma$, we call $\underline{n}_{K}(p)$ the canonical Weierstrass gap sequence (or simply, the Weierstrass gap sequence) of $p$. Accordingly, a $K$-Weierstrass point is called a canonical Weierstrass point, or simply, a Weierstrass point. 

The above definitions 
parallel their classical counterparts on algebraic curves. Let $\mathcal{C}$ be a smooth projective curve of genus $g \geq 2$ defined over a field of characteristic zero, and let $\mathcal{D}$ be a divisor of degree $d$ and rank~$r$ on $\mathcal{C}$. Then the locus $\WP_{\mathcal{C}}(\mathcal{D})$ of $\mathcal{D}$-Weierstrass points is finite, and the total weight satisfies the formula $\sum_{x \in  \WP_{\mathcal{C}}(\mathcal{D})} \wt_{\mathcal{D}}(x) 
= (r+1) (rg - r + d)$ (see~\cite[Theorem~6]{Lak81}). 

Motivated by a question of Baker~\cite[Remark~4.14]{Baker}, we seek to understand to what extent the classical theory of Weierstrass gap sequences can be extended to the tropical setting.

Our first result provides a characterization of the $D$-Weierstrass gap sequence at a point $p$ that is isolated in the $D$-Weierstrass locus. Consequently, the weight of such a point coincides with the definition used by Amini--Gierczak--Richman \cite{AGR}, and, together with their results, this yields a formula for the total weight of all $D$-Weierstrass points when the $D$-Weierstrass locus is discrete. For a point $p\in\Gamma$, denote by $D_p$ the $p$-reduced divisor in the linear equivalence class of $D$; see Section~\ref{sec:preliminaries}. The coefficient of $D_p$ at $p$ is denoted by $D_p(p)$.

\begin{Theorem}
\label{thm:Wloc:finite}
Let $\Gamma$ be a tropical curve of genus $g \geq 2$, and let 
$D$ be a divisor on $\Gamma$ of degree $d$ and nonnegative rank $r$. 
Then, for an isolated point $p$ in $\WL(D)$, 
the Weierstrass gap sequence at~$p$ with respect to $D$ is given by 
\[
\underline{n}_D(p)=(1, \, 2,\, \ldots,\, r, \, D_p(p) +1). 
\]
In particular, the $D$-Weierstrass weight of $p$ is $D_p(p)-r$. Furthermore, if $\WL(D)$ is finite, then 
\begin{equation}
\label{eqn:thm:Wloc:finite}
\sum_{p\in \WL(D)}
 \wt_{D}(p) = rg - r + d. 
\end{equation}
\end{Theorem} 

Since $\rank(K) = g-1$ and $\deg(K) = 2g-2$, 
Theorem~\ref{thm:Wloc:finite} implies that, if $\WL(K)$ is finite, then 
\begin{equation}
\label{eqn:thm:Wloc:finite:canonical}
\sum_{p\in \WL(K)}
 \wt(p) = g^2-1. 
\end{equation}

As an immediate corollary, we obtain an upper bound on the weight of the $D$-Weierstrass gap sequence of an isolated point. 
\begin{Corollary}
\label{cor:Wloc:finite}
Let $\Gamma$ be a tropical curve of genus $g \geq 2$, and let 
$D$ be a divisor on $\Gamma$ of nonnegative rank. Suppose that $p$ is an isolated point in $\WL(D)$. 
\begin{enumerate}
\item
If $D = K$, then $\wt(p) \leq g-1$.
\item
If $\deg D > 2g-2$, then $\wt_{D}(p) \leq g$. 
\item 
The bounds in \textup{(1)} and \textup{(2)} are optimal in the following sense: for \textup{(1)},  there exist a tropical curve $\Gamma$ of genus $g$ and a point $p \in \Gamma$ such that equality holds; for \textup{(2)}, for each $d$ with $d > 2g-2$, there exist a tropical curve $\Gamma$ of genus $g$, a divisor $D \in \Div(\Gamma)$ of degree $d$, and a point $p \in \Gamma$ such that equality holds.  
\end{enumerate}
\end{Corollary}

The situation becomes more subtle when $\WL(D)$ has a non-isolated connected component. The following proposition provides a sharp contrast to Corollary~\ref{cor:Wloc:finite}. 

\begin{Proposition}
\label{prop:contrast}
Let $\Gamma$ be a tropical curve of genus $g \geq 2$, let  
$D$ be a divisor on $\Gamma$ of nonnegative rank, 
and $p \in \Gamma$. 
\begin{enumerate}
\item
If $D = K$, then $\wt(p) \leq \frac{g(g-1)}{2}$. 
\item
If $\deg D > 2g-2$, then $\wt_{D}(p) \leq \frac{g(g+1)}{2}$. 
\item
The bounds in \textup{(1)} and \textup{(2)} are optimal, in the same sense as in Corollary~\ref{cor:Wloc:finite}. 
\end{enumerate}
\end{Proposition}

We explore the interplay between Weierstrass gap sequences on tropical curves and those on algebraic curves under tropicalization. 
Let $k$ be an algebraically closed, complete non-Archimedean valued field of characteristic zero, and let $\mathcal{C}$ be an algebraic curve of genus $g$ over $k$. Assume that 
$\mathcal{C}$ is a Mumford curve, meaning that it admits a strictly semi-stable model over the valuation ring of $k$ such that each irreducible component of the special fiber in this model is a smooth rational curve. Let $\Gamma$ denote the dual metric graph of the special fiber, and let $\tau\colon \mathcal{C}(k) \to \Gamma$ be the specialization map (see Section~\ref{subsec:Baker:specialization}). 

For each point $p\in \Gamma$, let $f_p$ be a rational function on $\Gamma$ which gives the reduced divisor $D_p$, that is, $D+\zero(f_p)=D_p$. For each tangent direction $\nu$ to $\Gamma$ at $p$, denote by $\slope_{p, \nu}(f_p) \in \ZZ$  the slope of $f_p$ at $p$ along $\nu$.

Let $B$ be a finite union of closed intervals of $\Gamma$, not necessarily a connected component of $\WL(D)$. We define the $D$-Weierstrass weight of $B$ analogous to the formula given in \cite{AGR} for the $D$-Weierstrass weight of a connected component of $\WL(D)$: 
\begin{equation}
\label{eqn:AGR:0:B}
\mu_D(B)
\colonequals \sum_{p \in B} D(p) +  (g(B)-c(B))\, r  
- \sum_{p \in \partial B} \sum_{\nu \in T_p\Gamma \cap \partial^{\rm out} B} \slope_{p, \nu}(f_p). 
\end{equation}
Here $D(p)$ denotes the coefficient of $D$ at $p$, $g(B)$ is the genus of $B$ (i.e., its first Betti number), $c(B)$ is the number of connected components of $B$, $\partial B$ is the boundary of $B$, and 
$\nu \in T_p\Gamma \cap \partial^{\rm out} B$ is a tangent direction that goes out of $B$ at $p$.

\begin{Theorem}
\label{thm:M:m}
Let $\Gamma$ be a tropical curve of genus $g \geq 2$ and let $D$ be a divisor on $\Gamma$ of nonnegative rank. 
Let $\tau\colon \mathcal{C}(k) \to \Gamma$ be as above and  
$\mathcal{D}$ be a divisor on $\mathcal{C}$ with 
$\tau(\mathcal{D}) = D$. 
Let $B$ be a finite union of closed intervals in $\Gamma$. 
Then we have 
\begin{equation}
\label{eqn:M:m}
\sum_{
\substack{
x \in {\rm WP}_{\mathcal{C}}(\mathcal{D}), 
\\
\tau(x) \in B
}
} \wt_{\mathcal{D}}(x)
\leq (r+1)\, \mu_D(B). 
\end{equation}
\end{Theorem}
When $B$ is a connected component of $\WL(D)$, it is shown in \cite{AGR} that equality holds. Theorem~\ref{thm:M:m} shows that, for a more general subset $B$, one obtains an inequality instead. Degeneration of Weierstrass points have been studied in \cite{Dia, EH, EM02, Baker, Amini, Gen21, AGR} among others. Theorem~\ref{thm:M:m} sheds new light on this phenomenon and yields interesting examples concerning degeneration of Weierstrass points on algebraic curves; see Section~\ref{sec:from:alg:to:trop}. 
We also get the following. 

\begin{Corollary}
\label{cor:M:m}
Notation as in Theorem~\ref{thm:M:m}, assume that $B = \{b\}$ is a singleton, $b \in \Gamma$.  Then $\mu_D(\{b\}) \leq \wt_D(b)$. In particular,
\[
\sum_{
x \in {\rm WP}_{\mathcal{C}}(\mathcal{D}), \; 
\tau(x) = b} \wt_{\mathcal{D}}(x)
\leq (r+1)\, \wt_D(b). 
\]
\end{Corollary}

For nonnegative integers $d$ and $r$, define
\begin{equation}
\label{eqn:def:S:g}
\mathcal{S}(r, d) \colonequals \bigl\{\,\underline{n} =  (n_1, \ldots, n_{r+1}) \in \ZZ^{r+1} \mid 1 \leq n_1 < \cdots < n_{r+1} \leq d+1\,\bigr\}. 
\end{equation}
We endow $\mathcal{S}(r, d)$ with the coordinatewise partial order $\leq$ (see Section~\ref{subsec:lex:coord} for details).  
We set $\WL_{\geq \underline{n}}(D)
\colonequals \bigl\{\,
p \in \Gamma \mid \underline{n}_D(p) \geq \underline{n}
\bigr\}$.
We prove the following (see Theorem~\ref{thm:upper}). 

\begin{Theorem}
\label{thm:intro:upper}
For each $\underline{n} \in \mathcal{S}(r, d)$, the set
$\WL_{\geq \underline{n}}(D)$ 
is closed and has finitely many connected components. 
\end{Theorem}

In connection with tropicalization, we note the following.

\begin{Proposition}
\label{prop:specialization:coordinatewise}
Let the notation be as in Theorem~\ref{thm:M:m}. 
Take any $x \in \mathcal{C}(k)$, and denote by $\underline{n}_{\mathcal{D}}(x)$ the Weierstrass gap sequence at $x$ with respect to $\mathcal{D}$ on $\mathcal{C}$. 
Then $\underline{n}_{D}(\tau(x)) \geq \underline{n}_{\mathcal{D}}(x)$.
\end{Proposition}

We say that a closed subset $B$ of $\Gamma$ is a {\em maximal $D$-Weierstrass locus} if there exists an $\underline{n} \in \mathcal{S}(r, d)$ such that $B$ is a connected component of $\WL_{\geq \underline{n}}(D)$ and $\underline{n}_{D}(p) = \underline{n}$ for all $p \in B$. 
By 
Theorem~\ref{thm:intro:upper}, there are only finitely many of such loci.

\smallskip

We now specialize to the case of canonical divisors. 
Since $\rank(K) = g-1$ and $\deg(K) = 2g-2$, for $g \geq 2$ we set
\begin{equation}
\label{eqn:def:S:g:canonical}
\mathcal{S}(g) \colonequals 
\mathcal{S}(g-1, 2g-2) =  
\bigl\{\left. \,\underline{n} =  (n_1, \ldots, n_{g}) \in \ZZ^{g} \;\right|\; 1 \leq n_1 < \cdots < n_{g} \leq 2g-1\,\bigr\}.  
\end{equation}
For any $p \in \Gamma$, 
the canonical Weierstrass gap sequence $\underline{n}_{K}(p) = (n_1, \ldots, n_g)$ belongs 
to~$\mathcal{S}(g)$. 
We call a maximal $K$-Weierstrass locus 
{\em a maximal Weierstrass locus}, and we call 
a point belonging to some maximal Weierstrass locus 
a {\em maximal Weierstrass point}. In particular, if $p$ is isolated in the Weierstrass locus, then $p$ is a maximal Weierstrass point. 

Motivated by Equation~\eqref{eqn:thm:Wloc:finite:canonical} in the case where the Weierstrass locus is finite, 
we consider the total weight of {\em maximal} Weierstrass points when $g = 2, 3$. We state here the genus $3$ case (see Section~\ref{subsec:genus:2} for the genus $2$ case.)  
We say that a tropical curve $\Gamma$ is {\em bridgeless} if removing any segment of $\Gamma$ consisting entirely of points of valence two yields a connected topological space. If, instead, such a segment $e$ exists in $\Gamma$ whose removal makes $\Gamma$ disconnected, then each two points $p, q$ on $e$ are linearly equivalent, so $\underline{n}_{K}(p)$ is constant on $e$. 
Therefore, in studying $\underline{n}_K(p)$ for $p \in \Gamma$, we may contract all such segments and assume that $\Gamma$ is bridgeless.

A point $v \in \Gamma$ is a {\em cut vertex} if  the deletion of $v$ from $\Gamma$ makes $\Gamma$ disconnected. 

\begin{Theorem}
\label{thm:max:WP:genus:3}
Let $\Gamma$ be a bridgeless tropical curve of genus $g = 3$ without cut vertices. Then each {\em maximal} Weierstrass locus is isolated. 
Let $\{p_1, \ldots, p_m\}$ denote the set of {\em maximal} Weierstrass points on $\Gamma$. 
Then we have the following. 
\begin{enumerate}
\item
Suppose that $\Gamma$ is hyperelliptic. Then $m = 4$ and 
we have 
\begin{equation}
\label{eqn:prop:max:WP:genus:3:a}
\sum_{i=1}^4
 \wt_{K}(p_i) = 12 = (g^3 - g)/2. 
\end{equation}
\item
Suppose that $\Gamma$ is non-hyperelliptic. 
 Then $4 \leq m \leq 8$,  and 
we have 
\begin{equation}
\label{eqn:prop:max:WP:genus:3:b}
\sum_{i=1}^m 
 \wt_{K}(p_i) = 8 = g^2-1. 
\end{equation}
\end{enumerate}
\end{Theorem}
See Section~\ref{sec:proof:thm:max:WP:genus:3} for the proof and for the treatment of the case where cut vertices are allowed. In particular, for a bridgeless tropical curve of genus $3$, whether or not cut vertices are present, maximal Weierstrass points are always finite in number, and the total sum of their weights assumes only very specific values. It is natural to ask whether this phenomenon persists in higher genus, that is, whether an analogous statement holds more generally. See also Question~\ref{q:max:W:1}.

\medskip
Finally, we consider classification of Weierstrass gap sequences. 
If $\mathcal{\mathcal{C}}$ is a smooth complex projective curve of genus $g \geq 2$ and $x \in \mathcal{C}(\CC)$ is a point, then the complement $\underline{n}_\mathcal{C}(x)^c$ in $\ZZ_{> 0}$ of the Weierstrass gap sequence $\underline{n}_\mathcal{C}(x)$ for $x$ forms a numerical semigroup of genus $g$, i.e.,  
$\left|\underline{n}_\mathcal{C}(x)\right| = g$ and $a + b \in \underline{n}_\mathcal{C}(x)^c$ for all $a, b \in \underline{n}_\mathcal{C}(x)^c$. Studies on which sequence is of the form $\underline{n}_\mathcal{C}(x)$ have a long history since Hurwitz \cite{Hu}. For example, Buchweitz (see \cite[p.~499]{EH}) showed that a certain numerical semigroup of $g = 16$ is not of the form $\underline{n}_\mathcal{C}(x)^c$ for any smooth complex projective curve $\mathcal{C}$ of genus $16$ and $x \in \mathcal{C}(\CC)$.  On the other hand, for $g \leq 8$, any numerical semigroup of genus $g$ is of the form $\underline{n}_\mathcal{C}(x)^c$ for some smooth complex projective curve $\mathcal{C}$ of genus $g$ and $x \in \mathcal{C}(\CC)$; cf. \cite{Komeda}. 

Motivated by the case of algebraic curves, we are interested in what sequence is of the form $\underline{n}_{K}(p)$ for a point $p$ on a tropical curve $\Gamma$. We set 
\begin{align*}
\mathcal{G}(g) & \colonequals  \left\{\underline{n}\in\mathcal{S}(g)  \;\left|\; 
\text{
$\underline{n} = \underline{n}_{K}(p)$ for some tropical curve $\Gamma$ of genus $g$ and $p \in \Gamma$}
\right.  \right\}, 
 \\
\mathcal{N}(g) & \colonequals  \left\{
\underline{n}\in\mathcal{S}(g)  \;\left|\;
\text{The complement $\underline{n}^c$ in $\ZZ_{>0}$ is a numerical semigroup} 
\right. \right\}, 
\end{align*}
where $\mathcal{S}(g)$ is defined in \eqref{eqn:def:S:g:canonical}.
(Here, the symbol $\mathcal{G}$ stands for ``gap,'' and 
$\mathcal{N}$ for ``numerical.'') 

\begin{Theorem}
\label{thm:main:1}
For $2 \leq g \leq 4$, $\mathcal{G}(g)$ is given as follows. 
\begin{enumerate}
\item
$\mathcal{G}(2) = \mathcal{N}(2) = \left\{(1, 2), (1, 3)\right\}$. 
\item
$\mathcal{G}(3) = \mathcal{N}(3) = \left\{(1, 2, 3), (1, 2, 4), (1, 2, 5), (1, 3, 5)\right\}$. 
\item
$\mathcal{G}(4) = \mathcal{N}(4) \cup \left\{(1, 2, 4, 6), (1, 2, 5, 6)\right\}$. 
Namely, we have 
\begin{align*}
\mathcal{G}(4)  & = 
\left\{
(1, 2, 3, 4), (1, 2, 3, 5), (1, 2, 3, 6), (1, 2, 3, 7), \right. 
\\
& 
\qquad 
\left. (1, 2, 4, 5), (1, 2, 4, 6), (1, 2, 4, 7),  (1, 2, 5, 6), (1, 3, 5, 7)  \right\}. 
\end{align*}
\end{enumerate}
\end{Theorem}

In particular, unlike the situation for algebraic curves, the complement $\underline{n}_{K}(p)^c$ of $\underline{n}_{K}(p)$ in $\ZZ_{> 0}$ is not necessarily a numerical semigroup. Similar tropical hyperelliptic curves that we use for the proof of 
Theorem~\ref{thm:main:1} yield such canonical Weierstrass sequences for any $g \geq 4$; see Remark~\ref{rmk:no:numerical:semigp}. 

We note that in \cite{Bo}, that we became aware of after finishing the writing of this paper,
Borz\`i considers  (canonical) Weierstrass sets on finite graphs and provides a sufficient condition that a sequence appears as a canonical Weierstrass gap sequence, see \cite[Thm.~C]{Bo}. In particular, his result yields that $\mathcal{G}(g) \setminus \mathcal{N}(g)$ is nonempty in general. Therefore, there is some overlap between Theorem~\ref{thm:main:1} and \cite[Thm.~C]{Bo}. 

\subsection*{Organization of this paper}
In the preliminary Section~\ref{sec:preliminaries}, we recall notions that will be used in the subsequent sections. 
This includes a brief review of the theory of divisors on tropical curves and tropical Abel--Jacobi maps, tropical Clifford's inequality, computation of ranks of divisors on hyperelliptic curves, weights of Weierstrass loci by Amini--Gierczak--Richman, Baker's specialization lemma, and specialization of Weierstrass points. 
We then study in Section~\ref{sec:Sg:etc} basic properties of 
tropical gap sequences, and prove Theorem~\ref{thm:intro:upper} and Proposition~\ref{prop:specialization:coordinatewise}. In 
Section~\ref{sec:proofs}, we give the proofs of Theorem~\ref{thm:Wloc:finite}, Corollary~\ref{cor:Wloc:finite}, and Proposition~\ref{prop:contrast}.  In Section~\ref{sec:from:alg:to:trop}, we study specialization from algebraic curves to tropical curves, and prove Theorem~\ref{thm:M:m} and 
Corollary~\ref{cor:M:m}. In Section~\ref{sec:proof:thm:max:WP:genus:3}, we prove Theorem~\ref{thm:max:WP:genus:3}, while in Section~\ref{sec:classification}, we prove 
Theorem~\ref{thm:main:1}. Finally, in Section~\ref{sec:open:q}, we list some open questions.

\bigskip
\noindent
{\sl Acknowledgment.}\quad
OA is partially supported by the project AdAnAr of the Agence Nationale de la Recherche (ANR). SK is partially supported by KAKENHI 23K03041 and 25H00587. 

%
\section{Preliminaries}
\label{sec:preliminaries}
In this section, we briefly review the theory of divisors on tropical curves. For details, see e.g. \cite{BJ} and references therein. 

\subsection{Divisors on tropical curves}
For a positive integer $m$ and a positive number $\varepsilon > 0$, let $S_m(\varepsilon)$ denote 
the star of radius $\varepsilon$ with $m$ branches. 
By a {\em tropical curve}, we mean a compact, connected metric space $\Gamma$ such that,  
for any $p \in \Gamma$, there are an open neighborhood $U$ of $p$ and an isometry 
$U$ with $S_m(\varepsilon)$ for some $m$ and $\varepsilon$ that sends $p$ to the center of $S_m(\varepsilon)$. We also call a tropical curve a metric graph, and we use both terms interchangeably. 
Further, we assume that  the metric on $\Gamma$ is given by 
the path metric, i.e., for points $p, q$ in $\Gamma$, the distance between $p$ and $q$ 
is the minimum length of any path from $p$ to $q$. The integer $m$ above 
is called the {\em valence} of $p$ and is denoted by $\val(p)$. The tangent space $T_p(\Gamma)$ is 
defined as the set of all outgoing tangent directions to $\Gamma$ at $p$, so $|T_p(\Gamma)| = \val(p)$. 
We usually identify a tangent direction $\nu$ at $p$ with the corresponding unit tangent vector to $\Gamma$ at $p$. 
The {\em genus} $g(\Gamma)$ of $\Gamma$ is 
defined as the first Betti number, i.e., $g(\Gamma) \colonequals \dim_\RR H_1(\Gamma, \RR)$. 

Suppose that $G$ is a loopless finite graph with the set of vertices $V(G)$ and the set of edges $E(G)$, 
and that $\ell\colon E(G) \to \RR_{> 0}$ is a length function. Then we obtain a tropical curve $\Gamma$ 
by associating to each $e \in E(G)$ a copy of the interval $[0, \ell(e)]$, with $\{0, \ell(e)\}$ identified with 
the endpoints of $e$, and then further identifying the ends of different intervals corresponding to 
the same vertex $v \in V(G)$. In this case, we say that $(G, \ell)$ is a loopless model of $\Gamma$. Each tropical curve in the sense of the previous paragraph arises in this way, admitting, in fact, infinity many models. Moreover, points of valence different from two belong to all models. 

Let $\Gamma$ be a tropical curve. 
Let $\Div(\Gamma)$ be the free $\ZZ$-module generated by points of $\Gamma$. 
An element of $\Div(\Gamma)$ is called a {\em divisor}, and is written as a finite sum $D = \sum_{p\in \Gamma} a_p (p)$ 
with $a_p \in \ZZ$ (with all but finitely many of the $a_p$ equal to zero). The coefficient $a_p$ at $p$ is also denoted by $D(p)$. 
The {\em support} of $D$ is defined by $\{p \in \Gamma \mid D(p) \neq 0\}$, denoted 
by $\Supp(D)$. 
The {\em degree} of $D$ is defined by 
$\deg(D) \colonequals \sum_{p \in \Gamma} D(p)$. The $\ZZ$-submodule of $\Div(\Gamma)$ consisting of degree $d$ divisors is denoted by $\Div^d(\Gamma)$. 
A divisor $D$ is {\em effective}, denoted by $D \geq 0$,  if $D(p) \geq 0$ for any $p\in \Gamma$. The {\em canonical divisor} $K$ of $\Gamma$ is defined by 
$K \colonequals \sum_{p \in \Gamma} (\val(p) - 2)(p) \in \Div(\Gamma)$.

A function $f\colon \Gamma \to \RR$ is called a {\em rational function} on $\Gamma$ if $f$ is piecewise affine with integer slopes. The space of rational functions on $\Gamma$ is denoted by $\Rat(\Gamma)$. 
(We note that it is natural to add $\infty$ to $\Rat(\Gamma)$ to make it into a tropical module. For simplicity, however, we do not include $\infty$ in this paper.)
For $f \in \Rat(\Gamma)$ and $p \in \Gamma$, 
we denote the slope of $f$ at the direction $\nu \in T_p(\Gamma)$ by 
$\slope_{p, \nu}(f)$ or simply by $\slope_\nu(f)$. 
Note that $\slope_{p, \nu}(f)$ is an integer and that 
except for a finitely many points $p \in \Gamma$, $\val(p) = 2$ 
and $\sum_{\nu \in T_p(\Gamma)} \slope_\nu(f) = 0$. 
The {\em principal divisor} of $f$ is defined by 
\begin{equation}
\label{eqn:zero:f}
\zero(f) \colonequals \sum_{p \in \Gamma} \left(-\sum_{\nu \in T_p(\Gamma)} \slope_\nu(f)\right) (p). 
\end{equation}

Two divisors $D, E \in \Div(\Gamma)$ are {\em linearly equivalent}, denoted by $D \sim E$, 
if $D - E = \zero(f)$ for some $f \in \Rat(\Gamma)$. 
For a divisor $D \in \Div(\Gamma)$, we set 
\begin{equation}
\label{eqn:Rat:D}
\Rat(D) \colonequals
\{f \in \Rat(\Gamma) \mid D + \zero(f) \geq 0\}
\end{equation}
and $|D| \colonequals \{E \in\Div(\Gamma) \mid \text{$D \sim E$ and $E \geq 0$}\}$.

The Baker--Norine {\em rank} (or simply the {\em rank}) of a divisor $D \in  \Div(\Gamma)$, denoted by $\rank(D)$, is defined as the maximum integer among $-1$ and integers $k \geq 0$  such that for all points $p_1, \ldots, p_k$ in $\Gamma$, the divisor $D - (p_1) -\cdots -(p_k)$ is linearly equivalent to an effective divisor. 
In fact, Luo \cite{Luo} shows that there exists a finite subset $W$ of $\Gamma$ such that 
the $\rank(D)$ is equal to the maximum integer among $-1$ and integers $k \geq 0$  such that for $p_1, \ldots, p_k \in W$, the divisor $D - (p_1) -\cdots -(p_k)$ is linearly equivalent to an effective divisor. Such a finite set $W$ is called a {\em rank-determining set} of $\Gamma$. 

With the Baker--Norine rank, we have the tropical Riemann--Roch formula. 

\begin{Theorem}[Tropical Riemann--Roch formula {\cite{BN, GK, MZ}}]
\label{thm:RR}
Let $\Gamma$ be a tropical curve of genus $g$. Then, for any $D \in \Div(\Gamma)$, we have 
\[
\rank(D) - \rank(K - D) = 1 - g + \deg(D).  
\]
\end{Theorem}

From the tropical Riemann--Roch formula, one can deduce 
Clifford's inequality. We say that a divisor $D \in \Div(\Gamma)$ is {\em special} 
if $\rank(D) \geq 0$ and $\rank(K-D) \geq 0$. 

\begin{Corollary}[Tropical Clifford inequality, {see \cite[Corollary 3.5]{BN}}]
\label{cor:RR}
Let $D$ be a special divisor. 
Then $\rank(D) \leq \frac{1}{2} \deg(D)$. 
\end{Corollary}

A tropical curve $\Gamma$ of genus $g \geq 2$ is {\em hyperelliptic} if there exists a divisor on $\Gamma$ of degree $2$ and rank $1$. 

The equality condition for the tropical Clifford inequality 
is a theorem of Coppens \cite{Coppens}. 

\begin{Theorem}[Coppens \cite{Coppens}]
\label{thm:Clifford}
Let $\Gamma$ be a tropical curve of genus $g \geq 4$, and let $D$ be a special divisor with $2 \leq \rank(D) \leq g-2$. If the equality $\rank(D) = \frac{1}{2} \deg(D)$
holds, then $\Gamma$ is hyperelliptic. 
\end{Theorem}

\subsection{Rank on hyperelliptic tropical curves}
\label{subsec:trop:hyper}
Let $\Gamma$ be a hyperelliptic tropical curve of genus $g \geq 2$. 
A {\em leaf edge} of $\Gamma$ is an edge $e$ of some loopless model $(G, \ell)$ of $\Gamma$ such that 
$e$ has an endpoint $p$ with $\val(p) = 1$. Let $\Gamma^\prime$ be the tropical curve without leaf edge 
obtained by successively contracting leaf edges of $\Gamma$. Then there exists an involution $\iota$ on $\Gamma^\prime$ such that $\Gamma^\prime/\langle \iota \rangle$ is a tree, and we take $v_0 \in \Gamma^\prime\; (\subseteq \Gamma)$ such that $\iota(v_0) = v_0$ (see e.g. \cite[\S3]{KY}). For an effective divisor $D$ on $\Gamma$, we set 
$
p_\Gamma(D) = \max\{k \in \ZZ_{\geq 0}\mid |D - 2k\, (v_0)| \neq \emptyset\}
$. 

We will use the following criterion for the computation of ranks 
on hyperelliptic tropical curves. 

\begin{Proposition}[{\cite[Theorem~1.14]{KY}}]
\label{prop:rank:hyp:ell:KY}
Let $\Gamma$ be a hyperelliptic tropical curve of genus $g \geq 2$ and let $D$ be an 
effective divisor on $\Gamma$. Then 
\[
\rank(D) 
= 
\begin{cases}
p_{\Gamma}(D) & \text{\textup{(}if $\deg D - p_{\Gamma}(D) \leq g$\textup{)}}, 
\\
\deg(D) - g & \text{\textup{(}if $\deg D - p_{\Gamma}(D) \geq g+1$\textup{)}}.
\end{cases}
\]
\end{Proposition}

\subsection{Total weight of the Weierstrass loci}
Let $\Gamma$ be a tropical curve of genus $g$. 
For any connected closed subset $A \subseteq \Gamma$, we denote by $\partial A$ the 
set of boundary points of $A$, and by $\partial^{\rm out} A$ the union of 
the tangent directions going out of $A$ at all $p \in \partial A$. The number of 
tangent directions going out of $A$ at $p$ is denoted by $\outdeg_A(p)$. In other words, 
$\outdeg_A(p) = \#(T_p(\Gamma) \cap \partial^{\rm out} A)$. 

Fix a point $q \in \Gamma$. A divisor $D \in \Div(\Gamma)$ is called $q$-reduced if it satisfies the following two conditions:
\begin{enumerate}
\item[(i)]
$D(p)\geq 0$ for all $p \in \Gamma\setminus \{q\}$. 
\item[(ii)]
For every connected closed subset $A \subseteq \Gamma\setminus \{q\}$, 
there exists $p \in \partial A$ such that $D(p) < \outdeg_A(p)$. 
\end{enumerate}

The following theorem is useful in computing the rank of divisors, see Section~5.1 of the survey paper~\cite{BJ}.

\begin{Theorem}\label{thm:reduced-non-negative-rank}
For each $D \in \Div(\Gamma)$, there exists a unique $q$-reduced divisor 
that is linearly equivalent to $D$. Further, $\rank(D) \geq 0$ if and only if the coefficient of the corresponding $q$-reduced divisor at $q$ is nonnegative. 
\end{Theorem}

We denote by $D_q$ the $q$-reduced divisor that is linearly equivalent to $D$. 

It follows that if $\rank(D)\geq 0$, then $D_q$ belongs to $|D|$. We denote $f_q$ a function in $\Rat(D)$ (see \eqref{eqn:Rat:D}) such that $D + \zero(f_q) = D_q$. Then for any tangent direction $\nu \in T_q(\Gamma)$ and any $f \in \Rat(D)$, we have 
\begin{equation}
\label{eqn:slope:fp}
\slope_{q, \nu}(f) \geq \slope_{q, \nu}(f_q). 
\end{equation}
Indeed, we set $g \colonequals f- f_q$. Then $D_q + \zero(g) \geq 0$. 
By \cite[Lemma~7]{Amini} and \cite[Lemma~A.5]{BakerShokrieh}, $g$ takes its minimum value at $q$. Thus $\slope_{q, \nu}(g) \geq 0$, from which the inequality \eqref{eqn:slope:fp} follows. (Note that in \emph{loc.~cit.} there is a sign difference in the definition of divisors of rational functions.) 

Let $D$ be a divisor on $\Gamma$ of 
degree $d$ and nonnegative rank $\rank(D) = r$. 
Let $A$ be a connected component of the Weierstrass locus $\WL(D)\colonequals \{p \in \Gamma \mid \rank(K - (r+1)(p)) \geq 0\}$ (See Section~\ref{subsec:D:Weierstrass} for details). Note that $A$ is a connected closed set of $\Gamma$, thus $A$ is a tropical curve. Let $g(A)$ denote the genus of $A$. For $p \in \partial A$, we take $f_p \in \Rat(\Gamma)$ such that 
$D + \zero(f_p) = D_p$. 
Let 
\begin{equation}
\label{eqn:AGR}
\mu_D(A) 
= \sum_{p \in A} D(p) +  (g(A)-1)\, r
+ \sum_{p \in \partial A} \sum_{\nu \in T_p(\Gamma) \cap \partial^{\rm out} A} \slope_{p, \nu}(f_p). 
\end{equation}

Amini, Gierczak and Richman \cite{AGR} show the following formula of the total weight of $\mu_D(A)$. 

\begin{Theorem}[Amini--Gierczak--Richman {\cite[Theorem~1.5]{AGR}}]
\label{thm:AGR}
Let $\Gamma$ be a tropical curve of genus $g \geq 2$, 
and let $D$ be a divisor on $\Gamma$ of 
degree $d$ and nonnegative rank $r$. 
Let $A _1\ldots, A_m$ be the set of all connected components of the Weierstrass locus $\WL(D)$. Then 
\[
\sum_{i=1}^m \mu_D(A_i) = r g - r + d.
\] 
\end{Theorem}

\subsection{Baker's specialization lemma and specialization of Weierstrass points}
\label{subsec:Baker:specialization}
Recall Baker's specialization lemma. 
Let $k$ be an algebraically closed, complete non-Archimedean  valued field, and $\mathcal{C}$ an algebraic curve of genus $g$ over $k$. Assume that 
$\mathcal{C}$ is a Mumford curve (see the paragraph after Proposition~\ref{prop:contrast} in Introduction). 
Let $\Gamma$ be the dual metric graph of the special fiber of the model, and let $\tau\colon \mathcal{C}(k) \to \Gamma$ be the specialization map. 

One way to see this specialization map is to use the analytification $\mathcal{C}^{\rm an}$ of $ \mathcal{C}$ in the sense of Berkovich. There $\Gamma$ is identified with the skeleton in $\mathcal{C}^{\rm an}$ of a strictly semi-stable model of $\mathcal{C}$ over the valuation ring of $k$, and there is a canonical retraction $\mathcal{C}^{\rm an} \to \Gamma$ taking a connected component of $\mathcal{C}^{\rm an} \setminus \Gamma$  to a unique boundary point in $\Gamma$. Noting that the set of $k$-valued points $\mathcal{C}(k)$ is included in $\mathcal{C}^{\rm an}$, the restriction of the retraction map to $\mathcal{C}(k)$ gives the specialization map $\tau$. 

Extending linearly, $\tau$ induces a map from divisors on $\mathcal{C}$ to divisors on $\Gamma$, which we still denote by $\tau\colon \Div(\mathcal{C}) \to \Div(\Gamma)$. 

\begin{Theorem}[Baker's specialization lemma {\cite[Lemma~2.8]{Baker}}]
\label{thm:Baker:sp}
For any divisor $\mathcal{D}$ on $\mathcal{C}$, 
we have $\rank(\tau(\mathcal{D})) \geq \rank_{\mathcal{C}}(\mathcal{D})$, 
where $\rank_{\mathcal{C}}(\mathcal{D}) = \dim_k H^0(\mathcal{C}, \OO_\mathcal{C}(\mathcal{D})) -1$. 
\end{Theorem}

In the rest of this subsection, we explain specialization of 
Weierstrass points, following the exposition in~\cite[\S5]{AGR}; see also \cite[\S12.2]{BJ}. We assume that $k$ is of characteristic zero. 

Let $\mathcal{D}$ be a divisor on $\mathcal{C}$ with rank $r$, where 
$r \colonequals \dim_k H^0(\mathcal{C}, \OO_\mathcal{C}(\mathcal{D})) -1$. 
Given a point $x \in \mathcal{C}(k)$, we define  
the {\em vanishing set} $S_x(\mathcal{D})$ of $\mathcal{D}$ at $x$ to be the set of orders of vanishing of global sections of 
$\OO_{\mathcal{C}}(\mathcal{D})$ at $x$. We have $|S_x(\mathcal{D})| = r+1$, and we write $\underline{n}_{\mathcal{D}}(x) = (n_1, \ldots, n_{r+1})$, where $n_1 < \cdots < n_{r+1}$ and $n_i - 1 \in 
S_x(\mathcal{D})$ for any $1 \leq i \leq r+1$. We call $\underline{n}_{\mathcal{D}}(x)$ the $\mathcal{D}$-gap sequence. 
The weight of the $\mathcal{D}$-gap sequence 
is defined as $\wt_{\mathcal{D}}(x) \colonequals \sum_{i=1}^{r+1} (n_i - i)$. 

A point $x \in \mathcal{C}(k)$ is a {\em $\mathcal{D}$-Weierstrass point} if $\underline{n}_{\mathcal{D}}(x) \neq (1, 2, \ldots,  r+1)$. We denote by $\WP_{\mathcal{C}}(\mathcal{D})$ 
the set of $\mathcal{D}$-Weierstrass points, which is a finite set. The {\em Weierstrass divisor} for $\mathcal{D}$ is defined by 
\begin{equation}
\label{eqn:mathcalW}
\mathcal{W} \colonequals \sum_{x \in \WP_{\mathcal{C}}(\mathcal{D})} \wt_{\mathcal{D}}(x) (x) \in \Div(\mathcal{C}). 
\end{equation}

Let $K_{\mathcal{C}}$ denote a canonical divisor on $\mathcal{C}$. 
If we fix a basis $\mathcal{F}$ of $H^0(\mathcal{C}, \OO_{\mathcal{C}}(\mathcal{D}))$, the corresponding Wronskian ${\rm Wr}_{\mathcal{F}}$ is a nonzero global section of 
$\OO_{\mathcal{C}}\left((r +1) \mathcal{D} + \frac{r(r+1)}{2} K_{\mathcal{C}}\right)$, 
whose divisor is precisely~$\mathcal{W}$. 

We seek an explicit formula for $\tau(\mathcal{W}) \in \Div(\Gamma)$. For this, we define 
$\Rat(\mathcal{D}) = \{F \in \Rat(\mathcal{C}) \mid \mathcal{D} + \zero(F) \geq 0\}$. For $p \in \Gamma$ and a tangent direction $\nu \in T_p(\Gamma)$, define $S^\nu(\mathcal{D})$ to be the set of integers occurring as $\slope_{p, \nu}({\rm trop}(F))$ for some $F \in \Rat(\mathcal{D})$. Here, identifying $\Gamma$ with a skeleton in $\mathcal{C}^{\rm an}$, 
${\rm trop}(F) \in \Rat(\Gamma)$ is defined to be $\rest{-\log |F|}{\Gamma}$, and satisfies $\tau(\mathcal{D}) + \zero({\rm trop}(F)) \geq 0$. We have $|S^\nu(\mathcal{D})| = r+1$, and we enumerate $S^\nu(\mathcal{D})$ as $s_0^\nu < \ldots < s_r^\nu$.  The following formula is proved in~\cite{Amini}. 

\begin{Theorem}[{\cite[Theorem~5.5]{AGR}}]
Notation as above, set $W \colonequals \tau(\mathcal{W})$ and $D \colonequals \tau(\mathcal{D})$ 
in $\Div(\Gamma)$. Then for any $p \in \Gamma$, we have 
\begin{equation}
\label{eqn:W(p)}
W(p) = (r+1)D(p)+\frac{r(r+1)}2K(p)-\sum_{\nu\in T_p\Gamma}\sum_{j=0}^r s^\nu_j.
\end{equation}
\end{Theorem}

\subsection{Effective locus in \texorpdfstring{$\Pic^d(X)$}{Pic(X)}}
\label{subsec:Pic:d}
Let $\Gamma$ be a tropical curve of genus $g \geq 1$. For each $d$, let $\Pic^d(\Gamma)$ be the linear equivalence classes of all degree $d$ divisors on $\Gamma$. For each divisor $D$ of degree $d$, denote by $[D]\in \Pic^d(\Gamma)$ its linear equivalence class. The effective locus $W_d\subseteq \Pic^d(\Gamma)$ is the locus of effective divisors classes, that is, 
\begin{equation}
\label{eqn:Wd}
W_d \colonequals \{[D] \in \Pic^d(\Gamma) \mid \rank(D) \geq 0\}. 
\end{equation}
The {\em tropical Jacobian} of $\Gamma$ is the  $g$-dimensional real torus 
\[
\Jac(\Gamma) \colonequals H_1(\Gamma, \RR)/H_1(\Gamma, \ZZ). 
\]
By tropical Abel--Jacobi theorem (see \cite{MZ} and \cite{BF}), there is an isomorphism
\[
\mathrm{AJ} \colon \Pic^0(\Gamma)\to \Jac(\Gamma).
\]

Fixing a point $p_0 \in \Gamma$, we get
a natural identification
\begin{align}\label{eqn:trop:AJ}
\Phi_{p_0}\colon \Pic^d(\Gamma) \to \Jac(\Gamma), \quad
\sum_{i=1}^k a_i (p_i) \mapsto \sum_{i=1}^k a_i  \mathrm{AJ}((p_i) -(p_0)). 
\end{align}
By an abuse of the notation, we define the {\em effective locus} $\widetilde{W}_d \subseteq \Jac(\Gamma)$ as the image of $W_d$ by $\Phi_{p_0}$. Changing the fixed point $p_0$ results in a translation of $\widetilde{W}_d$ in $\Jac(\Gamma)$.

The map $\Gamma \to \Jac(\Gamma)$ that sends $p$ to $\mathrm{AJ}((p)-(p_0))$ is piecewise affine, meaning that there exists a model of $\Gamma$ 
such that on each interval edge of this model, 
it factors through a piecewise affine map to the vector space 
$H_1(\Gamma, \RR)$. Since the composite $\Gamma^d  \to \Jac(\Gamma)$ is still piecewise affine, the effective locus $\widetilde{W}_d$ is a closed polyhedral subset of $\Jac(\Gamma)$.  It is proved  in \cite[Theorem~8.3]{GST22} that $\widetilde{W}_d$ is of pure dimension $d$ for each $0 \leq d \leq g$. 

We regard $\Pic^d(\Gamma)$ as a real torus via the bijection \eqref{eqn:trop:AJ}.  It thus follows that $W_d$ is a closed pure-dimensional polyhedral subset of $\Pic^d(\Gamma)$.

\section{\texorpdfstring{$D$}{D}-Weierstrass gap sequences and coordinatewise partial order}
\label{sec:Sg:etc}
In this section, we present basic properties of $D$-Weierstrass gap sequences. 

\subsection{\texorpdfstring{$D$}{D}-Weierstrass points and \texorpdfstring{$D$}{D}-Weierstrass gap sequences}
\label{subsec:D:Weierstrass}
Let $\Gamma$ be a tropical curve of genus $g \geq 0$. Let $D \in \Div(\Gamma)$. We set $d = \deg(D) \in \ZZ$ and $r = \rank(D) \geq -1$. 

\begin{Definition}
\label{def:D:Wei:gap}
For $p \in \Gamma$, we say that a positive integer $m$ is a 
{\em $D$-Weierstrass gap} at $p$ if 
$\rank(D - m(p)) < \rank(D - (m-1)(p))$. 
\end{Definition}

\begin{Lemma}
Let $p \in \Gamma$. Then there are exactly 
$r+1$ $D$-Weierstrass gaps at $p$, and the largest 
$D$-Weierstrass gap is at most $d+1$. 
\end{Lemma}

\Proof
If $r = -1$, then there is nothing to prove. Assume that $r \geq 0$. 
If $n \geq d+2$, then $\deg(D - (n-1)(p)) < 0$, so 
$\rank(D - n(p)) = \rank(D - (n-1)(p)) = -1$.  
This implies the second assertion. For the first assertion, 
we have 
\[
r = \rank(D) \geq \rank(D - (p)) \geq \rank(D - 2 (p))
\geq \cdots \geq \rank(D - (d+1)(p)) = -1. 
\]
Since removing a point from a divisor decreases the rank by at most one, there are exactly $r+1$ $D$-Weierstrass gaps. 
\QED

We denote the {\em $D$-Weierstrass gap sequence} of $p$ by 
$\underline{n} \colonequals (n_1, \ldots, n_{r+1})$ with $1 \leq n_1 < \cdots < n_{r+1} \leq d+1$, where $n_1, \ldots, n_{r+1}$ are the $D$-Weierstrass gaps at $p$. 
When we wish to emphasize $D$ and $p$, we write $\underline{n}_D(p)$ in place of $\underline{n}$. Note that from the definition, 
\begin{equation}
\label{eqn:r:nj:j}
\rank(D - n_i (p)) = r - i 
\qquad \text{for $1 \leq i \leq r+1$}. 
\end{equation}
The weight of the $D$-Weierstrass gap sequence $\underline{n}$ is defined by 
\[
\wt(\underline{n}) \colonequals \sum_{i=1}^{r+1} (n_i - i) \in \ZZ_{\geq 0}.
\] 

A point $p\in \Gamma$ is a {\em $D$-Weierstrass point} if 
$\underline{n}_D(p) \neq (1, 2, \ldots, r+1)$. This is equivalent to $\wt(\underline{n}(p))>0$. 

This notion of Weierstrass point is equivalent to the one in~\cite{Baker} and~\cite{AGR}, as we show next.

\begin{Lemma}
\label{lem:non-W:consecutive}
A point $p \in \Gamma$ is a $D$-Weierstrass point if and only if $D_p(p)\geq r+1$.  
\end{Lemma}

\Proof
If $\underline{n}_D(p) = (1, 2, \ldots, r+1)$, then we have $\rank(D - (r+1)(p)) = r - (r+1) = -1$ by \eqref{eqn:r:nj:j}, which by Theorem~\ref{thm:reduced-non-negative-rank} implies that $D_p(p) < r+1$.

Conversely, suppose that $p \in \Gamma$ verifies $D_p(p) < r+1$.  This means that $\rank(D - (r+1)(p)) = -1$. Since removing a point from a divisor decreases the rank by at most one, and 
\[
r = \rank(D) \geq \rank(D - (p)) \geq \rank(D - 2 (p))
\geq \cdots \geq \rank(D - (r+1)(p)) = -1, 
\]
we get $\underline{n}_D(p) = (1, 2, \ldots, r+1)$. 
\QED

Define the $D$-Weierstrass locus by 
\begin{equation}
\label{eqn:def:D:Wei:locus}
\WL(D)
\colonequals 
\{
p \in \Gamma \mid 
\underline{n}_D(p) \neq (1, 2, \ldots, r+1)
\}
= 
\{
p \in \Gamma \mid \rank(D - (r+1)(p)) \geq 0
\}. 
\end{equation}

We rephrase the above definitions in the case where $D$ is the canonical divisor $K$. 

\begin{Lemma}
Suppose that $D = K$. 
\begin{enumerate}
\item
A point $p \in \Gamma$ is a Weierstrass point if and only if $\rank(g (p)) \geq 1$.
\item
An integer $n$ is a Weierstrass gap if and only if 
$\rank((n-1) (p)) = \rank(n (p))$. 
\item
The Weierstrass locus is given by 
$\WL(K) = \{p \in \Gamma \mid 
\rank(g (p)) \geq 1\}$. 
\end{enumerate}
\end{Lemma}

\Proof
By the tropical Riemann--Roch formula (see Theorem~\ref{thm:RR}), we have $\rank(K) = g-1$. 
Again by the tropical Riemann--Roch formula, we have 
$\rank(g (p)) = 1 + \rank(K - g (p))$, which gives Assertions~(1) and (3). We also have 
\[
\rank(n (p)) = 1 - g + n + \rank(K - n (p)), \quad 
 \rank((n-1) (p)) = 1 - g + (n-1) + \rank(K - (n-1) (p)).  
\]
Thus $\rank((n-1) (p)) = \rank(n (p))$ if and only 
if $\rank(K - (n-1) (p)) > \rank(K - n (p))$, 
which gives Assertion~(2).  
\QED

\subsection{Coordinatewise partial order}
\label{subsec:lex:coord}
For $d \geq 0$ and $r \geq 0$, we set 
\[
\mathcal{S}(r, d) \colonequals \{\underline{n} =  (n_1, \ldots, n_{r+1}) \in \ZZ^g \mid 1 \leq n_1 < \cdots < n_{r+1} \leq d+1\}. 
\]
We endow a partial order $\leq$ with $\mathcal{S}(r, d)$. 
Let $\underline{n} = (n_1, n_2, \ldots, n_{r+1})$,  $\underline{n}^\prime = (n_1^\prime, n_2^\prime, \ldots, n_{r+1}^\prime) \in \mathcal{S}(r, d)$. We write 
\[
\underline{n} \leq \underline{n}^\prime
\quad\text{provided that}\quad
n_i \leq n_i^\prime \quad\text{for any}\quad i = 1, \ldots, r+1. 
\] 

A function $f\colon \Gamma \to \mathcal{S}(r, d)$ is called {\em upper-semicontinuous} 
if for each $\underline{n}\in \mathcal{S}(r, d)$, the subset $\{p \in \Gamma \mid f(p) \geq \underline{n}\}$ is closed in $\Gamma$.  

Let $\Gamma$ be a tropical curve and $D$ a divisor on $\Gamma$ of degree $d$ and rank $r$. 
For each $\underline{n} \in \mathcal{S}(r, d)$, we set 
\begin{align*}
& \WL_{\geq \underline{n}}(D)
\colonequals \bigl\{\,
p \in \Gamma \mid \underline{n}_D(p) \geq \underline{n}
\bigr\}. 
\end{align*}
We are going to show that $\WL_{\geq \underline{n}}(D)$ 
has a finitely many connected components and that 
the map $\varphi: \Gamma \to \mathcal{S}(r, d)$ assigning each $p \in \Gamma$ the $D$-Weierstrass gap sequence $\underline{n}_{D}(p)\in \mathcal{S}(r, d)$ is upper-semicontinuous with respect to the coordinatewise partial order. 

We need several lemmas. 
Recall that $W_{d-r} \colonequals \{[F] \in \Pic^{d-r}(\Gamma) \mid \rank(F) \geq 0\} \subseteq \Pic^{d-r}(\Gamma)$ is the {\em effective locus} (see Section~\ref{subsec:Pic:d}). We set    
\[
W_{r, d}
\colonequals \{[D] \in \Pic^d(\Gamma) 
\mid \rank(D) \geq r\} \subseteq \Pic^d(\Gamma). 
\]

\begin{Lemma}
\label{lem:Wrd:polyhedron}
The set $W_{r, d}$ is a closed polyhedral subset in $\Pic^d(\Gamma)$. 
\end{Lemma}

\Proof
Take a rank-determining set $W \subseteq \Gamma$ (see the paragraphs before 
Theorem~\ref{thm:RR}) and 
an effective divisor $E \in \Div(\Gamma)$ of degree $r$  with $\Supp(E) \subseteq W$. 

Let $\phi_E\colon \Pic^{d}(\Gamma) \to \Pic^{d-r}(\Gamma)$ be the isomorphism given by 
$D \mapsto D-E$. As we explained in Section~\ref{subsec:Pic:d}, $W_{d-r}$ is a closed polyhedral subset of $\Pic^{d-r}(\Gamma)$. 

Since $W$ is rank-determining, we have 
\[
W_{r, d} = \bigcap_{\substack{E\in \Div(\Gamma), 
E \geq 0
\\ 
\deg(E) = r, \, 
\Supp(E) \subseteq W}} \phi_E^{-1}(W_{d-r}). 
\]
Since the isomorphism $\phi_E^{-1}$ is a translation, $\phi_E^{-1}(W_{d-r})$ is a closed polyhedral subset of $\Pic^{d}(\Gamma)$, and we see that $W_{r, d}$ is a closed polyhedral subset of $\Pic^d(\Gamma)$. 
\QED

Let $I$ be a closed interval of $\Gamma$, which we identify with $I=[a,b] \subseteq \RR$. 
Let $D$ be a divisor on $\Gamma$. For nonnegative integers $m$ and $k$, we set  
\[
I(m, k) \colonequals \{p\in I \mid \rank(D-m(p)) =k\}. 
\]
By an open interval of $I$ we mean the intersection of an open interval in $\RR$ with $I$. 

\begin{Lemma}
\label{lem:3:6} 
Suppose that for all $p\in I$, $\rank(D-m(p)) \geq k$. 
Then $I(m, k)$ is a finite union of open intervals in $I$.
\end{Lemma}

\Proof
We define the map $\varphi\colon I \to \Pic^{d-m}(\Gamma)$ by $p \mapsto D-m(p)$. 
Then 
\[
I(m, k) = \varphi^{-1}(W_{k, d-m} \setminus W_{k+1, d-m}) = \varphi^{-1}(W_{k, d-m}) \setminus \varphi^{-1}(W_{k+1, d-m}).
\]
Since $\varphi$ is piecewise affine (see Section~\ref{subsec:Pic:d}) 
and $W_{k, d-m}$ and $W_{k+1, d-m}$ are closed polyhedral subsets of $\Pic^{d-m}(\Gamma)$ by Lemma~\ref{lem:Wrd:polyhedron}, we obtain the assertion. 
\QED

\begin{Theorem}
\label{thm:upper}
Let $D$ be a divisor on $\Gamma$ of degree $d$ and rank $r$. 
Then the map $\varphi: \Gamma \to \mathcal{S}(r, d)$ assigning each $p \in \Gamma$ the $D$-Weierstrass gap sequence $\underline{n}_{D}(p)\in \mathcal{S}(r, d)$ is upper-semicontinuous with respect to the coordinatewise partial order.  Furthermore, for any $\underline{n} \in \mathcal{S}(r, d)$, the set $\WL_{\geq \underline{n}}(D)$ is a finite union of closed intervals of $\Gamma$. 
\end{Theorem}

\Proof
Let $\underline{n} = (n_1, \ldots, n_{r+1})\in \mathcal{S}(r, d)$. 
We claim that 
$\{p \in \Gamma \mid \underline{n}_D(p) \geq \underline{n}\}$ is a finite union of closed intervals of $\Gamma$. This gives the upper-semicontinuity of $\varphi$ at the same time. 

To show the claim, 
we proceed by the induction on $|\underline{n}| \colonequals n_1 + \cdots + n_{r+1}$. For the base case $\underline{n} = (1, \ldots, r+1)$, we have $\{p \in \Gamma \mid 
\underline{n}_{D}(p) \geq (1, \ldots, r+1)
\} = \Gamma$. Since $\Gamma$ is written as a finite union of closed intervals of $\Gamma$, we are done. 

In general, given $\underline{n} \in \mathcal{S}(r, d)$, 
we set 
\[
M \colonequals 
\left\{\left.\underline{m} = (m_1, \ldots, m_{r+1}) \in   \mathcal{S}(r, d)
\;\right|\;  
\underline{n} \geq \underline{m}, \; |\underline{m}| = |\underline{n}|-1
\right\}. 
\]
For each $\underline{m} \in M$, there exists a unique $k = k(\underline{m})$ 
with $1 \leq k \leq r+1$ such that $m_k = n_k - 1$. 
We write $\underline{n}_D(p)= (n_1(p), \ldots, n_{r+1}(p))$. 
We have 
\[
\{p \in \Gamma \mid \text{$\underline{n}_D(p) \geq \underline{n}$}\}
= \bigcup_{\underline{m} \in M}
\left\{p \in \Gamma \mid \underline{n}_D(p) \geq \underline{m},\; 
n_{k}(p) \neq m_k\right\}. 
\]
By induction, $\{p \in \Gamma \mid \underline{n}_D(p) \geq \underline{m}\}$ 
is a finite union of closed intervals of $\Gamma$. 

Let $I$ be one of the closed intervals. Since $n_k(p) \geq m_k$ for any $p \in I$, we have $\rank(D - m_k(p)) \geq r-k$ by \eqref{eqn:r:nj:j}. Further, $n_{k}(p) =  m_k$ if and only if 
$\rank(D - m_k(p)) = r-k$ by \eqref{eqn:r:nj:j}. It follows from 
Lemma~\ref{lem:3:6} that $\{p \in I \mid n_{k}(p) = m_k\} = 
I(m_k, r-k)$ is a finite union of open intervals in $I$. 
Thus the set $\left\{p \in \Gamma \mid \underline{n}_D(p) \geq \underline{m},\; 
n_{k}(p) \neq m_k\right\}$ is a finite union of closed intervals of $\Gamma$, 
and we are done.
\QED

\begin{Remark}
\label{rmk:uniformity}
The proof of Theorem~\ref{thm:upper} actually shows that there exists a constant $C = C(g, d)$ depending only on $g, d$ such that for each integer $r$ and $\underline{n} \in \mathcal{S}(r, d)$, each tropical curve $\Gamma$ of genus $g$, and each divisor $D$ on $\Gamma$ of degree $d$, 
there exist at most $C$ connected components 
in $W_{\geq \underline{n}}(D)$. 
\end{Remark}

We end this section by giving the proof of  Proposition~\ref{prop:specialization:coordinatewise}. We follow the notation in Section~\ref{subsec:Baker:specialization}. 

\medskip
{\sl Proof of Proposition~\ref{prop:specialization:coordinatewise}.}\quad
Let the notation be as in Theorem~\ref{thm:M:m}. 
Then we need to show that $\underline{n}_{D}(\tau(x)) \geq \underline{n}_{\mathcal{D}}(x)$. 
It suffices to show that for any integer $j$, 
$\rank_\mathcal{C}(\mathcal{D} - j(x)) \leq 
\rank_\Gamma(D - j(\tau(x)))$. 
This is a consequence of Baker's specialization lemma (see Theorem~\ref{thm:Baker:sp}). 
\QED

\section{Proofs of Theorem~\ref{thm:Wloc:finite}, 
Corollary~\ref{cor:Wloc:finite}, and Proposition~\ref{prop:contrast}}
\label{sec:proofs}
In this section, we prove Theorem~\ref{thm:Wloc:finite}, 
Corollary~\ref{cor:Wloc:finite}, and Proposition~\ref{prop:contrast}. 

\subsection{Proof of Theorem~\ref{thm:Wloc:finite}}
\label{subsec:pf:thm:Wloc:finite} 
Let $\Gamma$ be a tropical curve of genus $g \geq 2$, 
and let $D$ be a divisor on $\Gamma$ of degree $d$ and 
nonnegative rank $r$. Note that since $r \geq 0$, we have $|D| \neq \emptyset$, so $d \geq 0$. 

Let $A_i$ be a connected component of the $D$-Weierstrass locus $\WL(D)$ of $\Gamma$ (see \eqref{eqn:def:D:Wei:locus}), and let $p$ be a boundary point of $A_i$. 
Since $A_i$ is closed (see Theorem~\ref{thm:upper}), we have $p \in A_i$. 

We write $\underline{n}_D(p) = (n_1, n_2, \ldots, n_{r+1}) \in \mathcal{S}(r, d)$. 

First, we show that $n_{r+1} = D_p(p)+1$. Indeed, we have
$\rank\left(D -  D_p(p)(p)\right)= 0$ and for any $j \geq 1$, $\rank\left(D -  (D_p(p) + j)(p)\right)= -1$. Thus $D_p(p)+1$ is the largest Weierstrass gap at $p$, and we get $n_{r+1} = D_p(p)+1$. 

Next we show that $\underline{n}_D(p) = (1, 2, \ldots, r, D_p(p)+1)$. To derive a contradiction, we assume that $\underline{n}_D(p) \neq (1, 2, \ldots, r, D_p(p)+1)$. 
Then $n_r \geq r+1$, so $\rank(D- r(p)) > \rank(D- n_r(p)) = 0$ (see \eqref{eqn:r:nj:j}). 
Note that any point $q \in \Gamma$ very near to $p$, 
we have $(r+1)(q) \sim r(p) + (q^\prime)$ for some $q^\prime \in \Gamma$. Then  
$\rank(D - (r+1)(q)) = \rank(D -r(p) - (q^\prime))\geq  \rank(D -r(p)) -1 \geq 0$. 
It follows that $q$ is a  Weierstrass point, which contradicts with $p$ being 
a boundary point of $A_i$. We obtain 
$\underline{n}_D(p) = (1, 2, \ldots, r, D_p(p)+1)$. 
In particular, $\wt_{D}(p) = D_p(p) - r$. 

Suppose now that the $D$-Weierstrass locus $A_i$ is a singleton $\{p_i\}$. We take $f_i \in \Rat(\Gamma)$ such that $D + \zero(f_{p_i}) = D_{p_i}$. By \eqref{eqn:zero:f},  we have $D(p_i) - \sum_{\nu \in T_{p_i}\Gamma} \slope_{p_i, \nu}(f_{p_i}) = 
D_{p_i}(p_i)$. 
It follows from \eqref{eqn:AGR} that 
\[
\mu_D(\{p_i\})
= D(p_i) - r 
- \sum_{\nu \in T_{p_i}\Gamma} \slope_{p_i, \nu}(f_{p_i}) = 
D_{p_i}(p_i) - r = \wt_D(p_i). 
\]

Suppose each $D$-Weierstrass locus $A_i$ is a singleton $\{p_i\}$. 
We write $\{p_1, \ldots, p_m\}$ for the set of $D$-Weierstrass points of $\Gamma$. By \cite{AGR} (see Theorem~\ref{thm:AGR}), we have 
$
\sum_{i=1}^m \mu_D(\{p_i\}) = rg-r+d$. 
Then we obtain $\sum_{i=1}^m \wt_D(p_i) = rg-r+d$. 
\QED

The proof of Theorem~\ref{thm:Wloc:finite} gives the following. 

\begin{Proposition}
Let $\Gamma$ be a tropical curve of genus $g \geq 2$, 
$D$ a divisor on $\Gamma$ of degree $d$ and nonnegative rank $r$.
Let $\underline{n} = (n_1, \ldots, n_{r+1})\in \mathcal{S}(r, d)$. 
Let $p_i$ ($i = 1, 2, \ldots$) and $p$ be points on $\Gamma$ such that $\lim_{i\to\infty} p_i = p$ and $\underline{n}_D(p_i) = \underline{n}$ for all $i$. 
We write 
$\underline{n}_{D}(p) = (n_1(p), \ldots, n_{r+1}(p))$, 
which satisfies $\underline{n}_{D}(p) \geq \underline{n}$ by 
Theorem~\ref{thm:upper}. 
Assume that $\underline{n}_{D}(p) \neq \underline{n}$ 
and take the minimal $1 \leq k \leq r+1$ with $n_k < n_k(p)$. 
Then $n_k + 1 < n_{k+1}$. 
\end{Proposition}

\Proof
To avoid notational conflict, in this proof, we denote the divisor 
with support $p \in \Gamma$ and coefficient $1$ by $[p]$, in place of $(p)$. 

For $\ell = 1, \ldots, r+1$,
we have $\rank(D - (n_\ell-1)[p_i]) = r- \ell +1$ and 
$\rank(D - n_\ell [p_i]) = r- \ell$ (see \eqref{eqn:r:nj:j}). Similarly, we have 
$\rank(D - (n_\ell(p)-1)[p]) = r- \ell +1$ and 
$\rank(D - n_\ell(p) [p]) = r- \ell$. 
Since $n_k < n_k(p)$, 
we have $\rank(D - n_k [p]) \geq \rank(D - (n_k(p) -1) [p])
= r-k+1$. For any point $q \in \Gamma$ very near to $p$, 
we have $(n_k+1)[q] \sim n_k [p] + [q^\prime]$ for some $q^\prime \in \Gamma$. It follows that $\rank(D - (n_k+1)[q]) = \rank(D - n_k [p] - [q^\prime]) \geq 
\rank(D - n_k [p]) - 1 \geq r-k$. In particular, for large $i$, 
we have $\rank(D - (n_k+1)[p_i]) \geq r-k$. 
By \eqref{eqn:r:nj:j}, we have $n_k + 1 < n_{k+1}$. 
\QED

\begin{Example}
Let $\Gamma$ be a tropical curve of genus $4$. Then a point $p \in \Gamma$ with $\underline{n}_{K}(p) = (1, 2, 5, 6)$ is {\em not} the limit of $p_i \in \Gamma$ with $\underline{n}_{K}(p_i) = (1, 2, 4, 5)$. 
\end{Example}

\subsection{Proof of Corollary~\ref{cor:Wloc:finite}}
\label{subsec:pf:cor:Wloc:finite} 
Before we start the proof of Corollary~\ref{cor:Wloc:finite}, we introduce 
two classes of tropical curves (i.e., metric graphs). 

\begin{Example}[dipole graph]
\label{eg:dipole}
A dipole graph (also called a ``banana'' graph) $B_g$ of genus $g \geq2$ consists of two vertices joined by $g+1$ edges, possibly of different lengths. Let $v, v^\prime$ denote the vertices of $B_g$, and let $e_0, \ldots, e_g$ denote the edges of $B_g$. Let $\ell_i \in \RR_{>0}$ be the length of $e_i$. In this case, we denote $B_g$ by $B_g(\ell_0, \ldots, \ell_g)$. 
\end{Example}

\begin{Example}[wheel graph]
\label{eg:wheel}
A wheel graph of genus $g$ consists of vertices $w$ and $v_1, \ldots, v_g$, 
where $w$ and $v_i$ are joined by one edge for all $1 \leq i \leq g$, 
$v_i$ and $v_{i+1}$ are joined by one edge for $i = 1, \ldots, g-1$, 
and $v_g$ and $v_1$ are joined by one edge. Lengths of edges are possibly different. 
\end{Example}

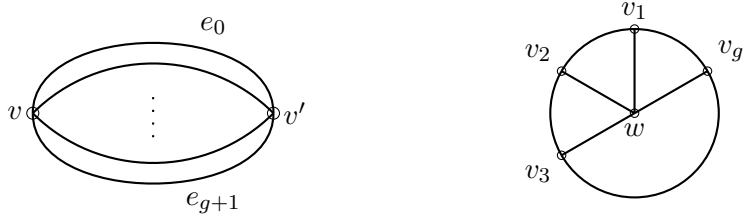
\begin{figure}[htb!]
\[
\begin{tikzpicture}[scale=0.8]
\draw[thick] (0,0) to [out=90,in=90,relative] (4, 0);
\draw[thick] (0,0) to [out=45,in=135,relative] (4, 0);
\draw[thick] (0,0) to [out=-45,in=-135,relative] (4, 0);
\draw[thick] (0,0) to [out=-90,in=-90,relative] (4, 0);
\draw[loosely dotted, thick] (2,-0.4)--(2, 0.4);
\draw (0,0)node[left]{$v$}; 
\draw (0, 0)circle (0.1);
\draw (4, 0)node[right]{$v^\prime$}; 
\draw (4, 0)circle (0.1);
\draw (3, 1.15)node[above]{$e_0$}; 
\draw (3, -1.15)node[below]{$e_{g+1}$}; 
\begin{scope}[shift={(10,0)}, scale=0.7]
\draw [thick] (0,0) circle [radius=2];
\draw [thick] (0,0)--(0, 2);
\draw [thick] (0,0)--(-1.732, 1);
\draw [thick] (0,0)--(-1.732, -1);
\draw [thick] (0,0)--(1.732, 1);
\draw (0, 0)circle (0.1);
\draw (0, 2)circle (0.1);
\draw (-1.732, 1)circle (0.1);
\draw (-1.732, -1)circle (0.1);
\draw (1.732, 1)circle (0.1);
\draw (0, 0)node[below]{$w$}; 
\draw (0, 2)node[above]{$v_1$}; 
\draw (-1.732, 1)node[above left]{$v_2$}; 
\draw (-1.732, -1)node[below left]{$v_3$}; 
\draw (1.732, 1)node[above right]{$v_g$}; 
\end{scope}
\end{tikzpicture}
\]
\caption{A dipole graph $B_g$ of genus $g$ (left) and a wheel graph 
of genus $g$ (right).}
\end{figure}

We begin the proof. 

(1) Recall that $\rank(K) = g-1 \;(\geq 1)$. 
Then Theorem~\ref{thm:Wloc:finite} gives 
\begin{equation}
\label{eqn:pf:cor:Wloc:finite}
\wt_{K}(p) = 
K_p(p) - \rank(K) \leq \deg(K) - (g-1) = g-1. 
\end{equation}

To show that this equality is optimal, consider 
the wheel graph of genus $g$, where each edge is assigned length one (see Example~\ref{eg:wheel}) and consider the central point $w$. Then $K \sim (2g-2) (w) = K_w(w)$. 

If $q \neq w \in \Gamma$ is very near to $w$, then 
$g(q) \sim (g-1)(w) + (q^\prime)$ for some $q^\prime \in \Gamma$ with $q^\prime \neq w$. 
Since $\outdeg_{\{w\}}(w) = g$, we see that 
 $(g-1)(w)  - (q^\prime)$ is a $q^\prime$-reduced divisor. Then 
\[
\rank(K - g (q)) = 
\rank((2g-2) (w) - (g-1)(w)  - (q^\prime)) = 
\rank((g-1)(w)  - (q^\prime)) = -1, 
\]
and $q$ is not a Weierstrass point. It follows that $w$ is an isolated 
Weierstrass point. 

Since $K_w(w) = 2g-2 = \deg(K)$, 
Equation~\eqref{eqn:pf:cor:Wloc:finite} is an equality. 

(2) Since $D$ has nonnegative rank by the assumption, we may assume that $D$ is 
effective. Since $\deg(D) > 2g-2$, the tropical Riemann--Roch formula (see Theorem~\ref{thm:RR}) implies that  
$\rank(D) = 1 - g + \deg(D) + \rank(K - D) = \deg(D)-g$. Then Theorem~\ref{thm:Wloc:finite} gives 
\begin{equation}
\label{eqn:pf:cor:Wloc:finite:2}
\wt_{D}(p) = 
D_p(p) - \rank(D) \leq \deg(D) - \rank(D) = g. 
\end{equation}

To show that this equality is optimal, consider 
the dipole graph $B_g$ of genus $g$, connecting two vertices $v$ and $v^\prime$ with $g+1$ edges (see Example~\ref{eg:dipole}). We assign each edge length one. For any integer $d > 2g-2$, we set $D \colonequals d(v) \in \Div(B_g)$. 

By the above computation, we have $\rank(D) = d-g$. 
Set $r \colonequals \rank(D)$. 
If $q \neq v \in B_g$ is very near to $v$, then 
$(r + 1)(q) \sim r (v) + (q^\prime)$ for some $q^\prime \in B_g$ with $q^\prime \neq v$. 
Since $\outdeg_{\{v\}}(v) = g+1$, we see that $g(v)  - (q^\prime)$ is a $q^\prime$-reduced divisor. Then 
\[
\rank(d (v) - (r + 1)(q)) = 
\rank(d (v) - r (v) - (q^\prime)) = 
\rank(g(v)  - (q^\prime)) = -1, 
\]
and $q$ is not a $D$-Weierstrass point. It follows that $v$ is an isolated 
$D$-Weierstrass point.

Since $D_v(v) = D(v) = d = \deg(D)$, 
Equation~\eqref{eqn:pf:cor:Wloc:finite:2} is an equality. 
\QED

\subsection{Proof of Proposition~\ref{prop:contrast}}
\label{subsec:pf:prop:contrast} 
We begin the proof of Proposition~\ref{prop:contrast}~(1)(3). 
Recall that $\rank(K) = g-1$. 
We write $\underline{n}_K(p) = (n_1, n_2, \ldots, n_{g})$. 
By \eqref{eqn:r:nj:j}, we have 
\[
\rank(K - (n_i-1) (p)) = (g-1) - (i-1), \quad
\rank(K -n_i (p)) = (g-1) -i 
\]
for all $i = 1, \ldots, g$. 
By tropical Clifford's inequality (see Theorem~\ref{thm:Clifford}), we get 
\[
\rank(K - (n_i-1) (p)) \leq \frac{1}{2}((2g-2) - (n_i-1)),
\]
so, $(g-1) - (i-1) \leq \frac{1}{2} ((2g-2) - (n_i-1))$, which implies that 
$n_i \leq 2i -1$. Then 
\[
\wt_{K}(p)
= \sum_{i=1}^g (n_i - i) \leq \sum_{i=1}^g (i-1) = \frac{g(g-1)}{2}. 
\]
To show that the bound is optimal, let $\Gamma$ be a hyperelliptic curve and let $p$ be a point with $\rank(2(p)) = 1$. 
Then $\underline{n}_{K}(p) = (1, 3, \ldots, 2g-1)$ and $\wt_{K}(p) =  \frac{g(g-1)}{2}$ (see Proposition~\ref{prop:rank:hyp:ell:KY}). 

\medskip
Next, we prove Proposition~\ref{prop:contrast}~(2)(3). 
 We set $r = \rank(D)$ and $d = \deg(D)$. 
We write $\underline{n}_D(p) = (n_1, n_2, \ldots, n_{r+1})$. 
Since $d > 2g-2$, the tropical Riemann--Roch formula (see Theorem~\ref{thm:RR}) gives $r = d-g$. The formula also implies that, 
for $0 \leq i \leq d-2g + 1$, 
$\rank(D - i(p)) = 1 - g + d - i + \rank(K - (D - i (p))) = d- g - i$, 
because $\deg(K - (D - i (p))) = 2g-2 - d + i \leq  2g-2 - d + (d-2g+1) = -1 < 0$. 
We set $\ell \colonequals d-2g + 1$. 

Then, for $1 \leq i \leq \ell$, $n_i = i$ and $\rank(D - i(p)) = r - i = d-g - i$. 
For $1 \leq j \leq g$, we define $m_j \geq 1 $ by $n_{\ell + j} = \ell + m_j$. Note that $\ell + g = d-g + 1 = r + 1$, so $n_{\ell + g} = n_{r+1}$. 
By the definition of the $D$-Weierstrass gap, we have 
\[
\rank(D - (\ell + m_j-1) (p)) = r - \ell -(j-1), \quad
\rank(D - (\ell + m_j) (p)) = r - \ell -j.  
\]

\begin{Claim}
We have $m_j \leq 2j$ for $1 \leq j  \leq g$. 
\end{Claim}

Indeed, if $m_j \leq j$, then trivially $m_j \leq 2j$. If $m_j \geq j+1$, then  
the tropical Riemann--Roch formula gives 
\begin{align*}
& \rank(K - (D - (\ell + m_j-1) (p))) \\
& \quad = \rank(D - (\ell + m_j-1) (p)) - (1 - g + \deg(D - (\ell + m_j-1) (p))) \\
& \quad = r-\ell - (j-1) - (1 - g +d - (\ell  +m_j-1))
= m_j - j - 1 \geq 0. 
\end{align*}
We also have $\rank(D - (\ell + m_j-1) (p)) = r - \ell -(j-1) 
= (d-g) - (d-2g + 1) - (j-1) = g-j \geq 0$. Thus $D - (\ell + m_j-1) (p)$ is special. 

By tropical Clifford's inequality (see Theorem~\ref{thm:Clifford}), we get 
\[
\rank(D - (\ell + m_j-1)(p)) \leq \frac{1}{2} \deg(D - (\ell + m_j-1)(p)),
\]
i.e., $g - j \leq \frac{1}{2} (d - (\ell + m_j-1))$, which implies that 
$m_j \leq 2j$. We have shown the claim. 

Then 
\[
\wt_{D}(p)
= \sum_{i=1}^\ell (n_i -i) + \sum_{j=1}^g (n_{\ell+j} -\ell-j)
= \sum_{j=1}^g (m_j-j) \leq \sum_{j=1}^g j 
= \frac{g(g+1)}{2}.  
\]

To show that the bound is optimal, let $\Gamma$ be a hyperelliptic curve and let $v_0$ be a point in the paragraph before Proposition~\ref{prop:rank:hyp:ell:KY}. For an integer $d > 2g-2$, we set $D = d (v_0) \in \Div(\Gamma)$. By Proposition~\ref{prop:rank:hyp:ell:KY}, we get 
\begin{align*}
& \rank(D - i(v_0))= d-g-i \; (1 \leq i \leq \ell), 
& \rank(D - (\ell  + j)(v_0))= \lfloor \frac{d-(\ell+j)}{2} \rfloor \; (1 \leq j \leq g). 
\end{align*}
We have $\underline{n}_D(p) = (1, 2, \ldots, \ell, \ell+2, \ell + 4, \ldots, \ell + 2g)$ and $\wt_{D}(p) =  \frac{g(g+1)}{2}$. 
\QED

\section{From algebraic curves to tropical curves}
\label{sec:from:alg:to:trop}
In this section, we prove Theorem~\ref{thm:M:m} and 
Corollary~\ref{cor:M:m}. 
We begin the proof of Theorem~\ref{thm:M:m}.

\bigskip
\noindent
{\sl Proof of Theorem~\ref{thm:M:m}.}\quad 
We follow the notation in Section~\ref{subsec:Baker:specialization}. 
In particular, $\WP_{\mathcal C}({\mathcal D})$ is the set of ${\mathcal D}$-Weierstrass points and ${\mathcal W}  \colonequals \sum_{x \in \WP_{\mathcal C}({\mathcal D})} \wt_{\mathcal D}(x) (x) \in \Div({\mathcal C})$ is the $\mathcal D$-Weierstrass divisor of $\mathcal D$ on $\mathcal C$, $W \colonequals \tau(\mathcal W)$, and $D \colonequals \tau(\mathcal D)$ in $\Div(\Gamma)$. Furthermore, given $p \in \Gamma$ and $\nu \in T_p(\Gamma)$, $s_0^\nu<\cdots<s_r^\nu$ are all the possible slopes at $p$ along $\nu$ of tropicalizations $\mathrm{trop}(F)$ for a rational function $F$ on $\mathcal C$ with ${\mathcal D}+\zero(F)\geq 0$. Since $\mathrm{trop}(F)$ satisfies
$D +\zero(\mathrm{trop}(F)) \geq 0$, by \eqref{eqn:slope:fp}, we have 
\begin{equation}
\label{eqn:s0}
s_j^\nu \geq s_0^\nu + j \quad (j = 0, 1, \ldots, r) 
\quad\text{and}\quad 
s_0^{\nu} \geq \slope_{p,\nu}(f_p). 
\end{equation}

Choose a loopless model  $G=(V,E)$ of $\Gamma$ such that $\Supp(W) \subseteq V$, and such that $B=\bigcup_{e\in F} e$ for the subset $F\subseteq E$ consisting of those edges in $E$ that connect two vertices in $V\cap B$ (in other words, $B$ is the induced subgraph on the vertices $V\cap B$). Without loss of generality, subdividing the edges of $G$ if necessary, we can assume that the sequence of slopes $s_0^\nu<\dots<s_r^\nu$ is constant over each edge $e\in E$, when $p \in e$ moves from an endpoint of $e$ to the other endpoint, and $\nu$ is the moving direction of $p$ on $e$. 

By Equations~\eqref{eqn:mathcalW} and~\eqref{eqn:W(p)}, taking the sum over $V \cap B$, we have  
\[
\sum_{
\substack{
x \in {\rm WP}_{\mathcal{C}}(\mathcal{D}), 
\\
\tau(x) \in B
}
} \wt_{\mathcal{D}}(x)
= 
\sum_{u\in V \cap B}	W(u) = (r+1)\sum_{u\in V \cap B}D(u)+\frac{r(r+1)}2\sum_{u\in V \cap B}K(u)-\sum_{u\in V \cap B}\sum_{\nu\in T_u\Gamma}\sum_{j=0}^r s^\nu_j.
\]
We have
\begin{align*}
	\sum_{u\in V\cap B}K(u) 
	&= \sum_{u\in V\cap B} (\val(u)-2) 
	\\
	& = 2|F|+|\{uw\in E \, |\, u\in B, w\notin B\}| - 2 |V\cap B| \\
	& \qquad\qquad \text{(since $uw \in F$ contributes two vertices $u, w \in V \cap B$)}
	\\
	&= 2(g(B) - c(B)) + |\{ uw\in E \, |\, u\in B, w\notin B\}|,
\end{align*}
where $c(B)$ is the number of connected components of $B$, the genus $g(B)$ equals $|F|-|V\cap B|+c(B)$, and we denote an edge $e \in E$ with endpoints $u, v$ by $uv$. 

With slight abuse of notation, we denote by $uv$ the tangent direction from $u$ to $v$ 
on an edge $uv$. Since by our assumption, the sequence of $(r+1)$ slopes along each edge of the model is constant, for each edge $uv$ in $F$, we have $s_j^{uv} = -s_{r-j}^{vu}$ for each $j=0, \dots, r$. This implies that
\[\sum_{j=0}^r s_j^{uv} 
+\sum_{j=0}^r s_j^{vu} =0.
\]
Therefore, 
\[
\sum_{u\in V\cap B}\sum_{\nu\in T_u\Gamma}\sum_{j=0}^r s^\nu_j = \sum_{uw\in E, u\in B, w\notin B} \sum_{j=0}^r s_j^{uw}.
\]
For each edge $uw$, using \eqref{eqn:s0}, we get
\[
	\sum_{j=0}^r s_j^{uw} \geq (r+1)s_0^{uw} +\frac{r(r+1)}2 \geq (r+1)\slope_{u, uw}(f_u) +\frac{r(r+1)}2.
\] 
Hence, we get
{\allowdisplaybreaks
\begin{align*}
	\sum_{u\in V\cap B}	W(u) 
	& \leq  (r+1)\sum_{u\in B}D(u) + r(r+1)(g(B) - c(B)) +\frac{r(r+1)}2|\{ uw\in E \, |\, u\in B, w\notin B\}| \\
	& - (r+1)\sum_{uw\in E, u\in B, w\notin B}\slope_{u, uw}(f_u) -\frac{r(r+1)}2|\{ uw\in E \, |\, u\in B, w\notin B\}|\\
	&= (r+1)\sum_{u\in B}D(u) + r(r+1)(g(B) - c(B)) - (r+1)\sum_{uw\in E, u\in B, w\notin B}\slope_{u, uw}(f_u)\\
	&=(r+1)\,\mu_D(B),
\end{align*}
}
see Equation~\eqref{eqn:AGR:0:B}. 
This completes the proof. 
\QED

Next, we prove Corollary~\ref{cor:M:m}. 

\bigskip
\noindent
{\sl Proof of Corollary~\ref{cor:M:m}.}\quad 
Let $f_b$ be an element of $\Rat(\Gamma)$ such that the divisor $D+\zero(f_b)$ is $b$-reduced. By the definition of $\zero(f_b)$ (see \eqref{eqn:zero:f}), we have 
$D_b(b) = D - \sum_{\nu \in T_b\Gamma} \slope_{\nu}(f_b)$. 
By the definition of $\mu_D(B)$ (see \eqref{eqn:AGR:0:B}), we have 
\[
\mu_D(\{b\}) = D(b) -  r  
- \sum_{\nu \in T_b\Gamma} \slope_{\nu}(f_b) = D_b(b)-r.  
\]
As we see in the proof of Theorem~\ref{thm:Wloc:finite}, in the gap sequence $\underline{n}_D(b)=(n_1, \ldots, n_{r+1})$, we have $n_{r+1} = D_b(b)+1$. 

For $1 \leq j \leq r$, we have $n_j \geq j$. Thus 
\[
\wt_D(b) = \sum_{j=1}^{r+1} (n_j - j) \geq n_{r+1} - (r+1) = D_b(b) - r = \mu_D(\{b\}), 
\]
as desired. The second assertion follows from 
the first assertion and Theorem~\ref{thm:M:m}. 
\QED

We illustrate Theorem~\ref{thm:M:m} with two examples. Note that, in the case of $D = K$, Theorem~\ref{thm:M:m} implies that 
$
\sum_{
x \in {\rm WP}_{\mathcal{C}}(K_{\mathcal{C}}), 
\tau(x) \in B
} \wt_{\mathcal{D}}(x)
\leq g \cdot \mu_{K}(B)$.

\begin{Example}
\label{eg:dipole:3}
Suppose that $\Gamma$ is a dipole graph of genus $3$, consisting of two vertices $v$ and $v'$ joined by $4$ edges $e_0, e_1, e_2, e_3$, possibly of different lengths (see Example~\ref{eg:dipole}). We identify $e_i$ with the closed interval $[0, \ell_i]$ with $v = 0$ and $v^\prime = \ell_i$ for $0 \leq i \leq 3$. 

We assume that the length of $e_0$ is $1$, so 
$e_0$ is identified with the closed interval $[0, 1]$. 
Here we consider the case where $D = K$. 
The Weierstrass gap sequences on $\Gamma$ are as in Table~\ref{table:B3}. 
Note that $\Gamma$ has no points with gap sequence $(1, 2, 5)$. 

\begin{table}[htb!]
\begin{tabular}{l|llllllllll}
\text{$\rank(k(p))$ for $1 \leq k \leq 5$} & $1$ & $2$ & $3$ & $4$ & $5$ &  & $\underline{n}_{K}(p)$ & & $\wt_{K}(p)$ 
\\ \hline
\smallskip
$p \in \left[0, \frac13\right) \cup \left(\frac23, 1\right]$  & $0$ & $0$ & $0$  & $1$ & $2$  & & $(1, 2, 3)$ && $0$\\ \hline
\smallskip
$p \in \left[\frac13, \frac23\right] \setminus \left\{\frac12\right\}$  & $0$ & $0$ & $1$ & $1$ & $2$  & & $(1, 2, 4)$ && $1$ 
\\  \hline
\smallskip
$p = \frac12$ & $0$ & $1$ & $1$ & $2$ & $2$  & & $(1, 3, 5)$ && $3$ \\
\end{tabular}
\caption{Weierstrass gap sequences on $B_3$}
\label{table:B3}
\end{table}

The connected component of the Weierstrass locus 
of $\Gamma$ on $e_0$ is equal to the closed interval $[-1/3, 1/3]$. 

Let $k$ be an algebraically closed, complete non-Archimedean valued field, and $\mathcal{C}$ a  Mumford curve of genus $3$ over $k$ with dual metric graph $\Gamma$. 
Let $\tau\colon \mathcal{C}(k) \to \Gamma$ be the specialization map.

\medskip
\emph{Case I: the curve is hyperelliptic}. \quad
The case where $\mathcal{C}$ is hyperelliptic  is clean. There are exactly $4(g-1) = 8$ Weierstrass points on $\mathcal{C}$, and each such point $x$ has the Weierstrass gap sequence $(1, 3, 5)$ and $\wt_\mathcal{C}(x) = \frac{g(g-1)}{2} = 3$.  

By Baker's specialization lemma (see Theorem~\ref{thm:Baker:sp}), 
$\tau(x)$ must be the middle point of each edge of $\Gamma$. 
So there are exactly $2$ Weierstrass points $x_1, x_2$ on $\mathcal{C}$ that land to the Weierstrass locus $A = \left[1/3, 2/3\right]\subseteq e_0$, and we have $\tau(x_1) = \tau(x_2) = \frac{1}{2} \in [0, 1] = e_0$. 

\medskip
\emph{Case II: the curve is non-hyperelliptic}. \quad
Now we assume that $\mathcal{C}$ is non-hyperelliptic. 

Let $x_1, \ldots, x_r \in \mathcal{C}$ be the Weierstrass points on $\mathcal{C}$ that land to $A = \left[1/3, 2/3\right]\subseteq e_0$.  By \eqref{eqn:AGR}, we compute  
$
\mu_{K}(A) = 2. 
$
By \cite{AGR} (see Theorem~\ref{thm:AGR}), we have $\wt_\mathcal{C}(x_1) + \cdots + \wt_\mathcal{C}(x_r) = g \cdot \mu_{K}(A) = 3 \cdot 2 = 6$. We arrange $x_1, \ldots, x_r$ so that $\tau(x_1) \leq \tau(x_2) \leq \cdots \leq \tau(x_r)$. 

We take $B = [1/3, 1/2-\varepsilon]$ for small $\varepsilon > 0$. 
By \eqref{eqn:AGR:0:B}, 
we compute 
$
\mu_{K}(B) = 1 
$. 
Theorem~\ref{thm:M:m} tells us 
that there are at most $3$ points $x_1, x_2, x_3$ that land to $[1/3, 1/2)$. 

For example, if we assume that there is no $x_i$ with $\tau(x_i) = 1/2$. If 
$\underline{n}_{K_\mathcal{C}}(x_i) = (1, 2, 5)$, then Proposition~\ref{prop:specialization:coordinatewise} implies that $\underline{n}_{K}(\tau(x_i)) \geq \underline{n}_{K_\mathcal{C}}(x_i)$ and since there is no point on $\Gamma$ with Weierstrass gap sequence $(1, 2, 5)$, we have $\underline{n}_{K}(\tau(x_i)) = (1, 3, 5)$ and $\tau(x_i) = 1/2$, which contradicts with our assumption. Thus $\underline{n}_{K_\mathcal{C}}(x_i) = (1, 2, 4)$ for all $1 \leq i \leq r$, and we get $r = 6$. By above, 
we see that $x_1, x_2, x_3$ land to $[1/3, 1/2)$ and 
$x_4, x_5, x_6$ land to $(1/2, 2/3]$. 
\end{Example}

\begin{Example}
\label{eg:hyperell:3:cut:vertices}
Let $\Gamma$ be a hyperelliptic curve of genus $3$ comprising of 
three circles $v, v^\prime, v^{\prime\prime}$, where the two vertices $v^\prime$ and  $v^{\prime\prime}$ are antipodal with respect to the middle circle. 
\begin{figure}[htb!]
\[
\begin{tikzpicture}[scale=0.7]
\draw[shift={(1,0)}] [thick] (0,0) circle [ radius=1];
\draw[shift={(1,0)}] [thick] (2,0) circle [ radius=1 ];
\draw[shift={(1,0)}] [thick] (4,0) circle [ radius=1 ];
\draw[shift={(1,0)}] (-1, 0)node[left]{$v$}; 
\draw[shift={(1,0)}] (-1, 0)circle (0.1);
\draw[shift={(1,0)}] (1, 0)node[right]{$v^\prime$}; 
\draw[shift={(1,0)}] (1, 0)circle (0.1);
\draw[shift={(1,0)}] (3, 0)node[right]{$v^{\prime\prime}$}; 
\draw[shift={(1,0)}] (3, 0)circle (0.1);
\end{tikzpicture}
\]
\caption{Chain of three circles.}
\end{figure}
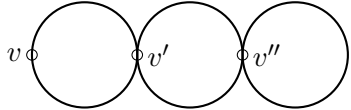

Let $e$ be a segment connecting 
$v$ and $v^\prime$. We assume that the length of $e$ is one, and we identify $e$ with 
the closed segment $[0, 1]$ with $v = 0$ and $v^\prime = 1$. 

Let $\mathcal{C}$ a  Mumford curve of genus $3$ with dual metric graph $\Gamma$. 
Let $\tau\colon \mathcal{C}(k) \to \Gamma$ be the specialization map.   

Since every point in $\Gamma$ is a Weierstrass point (cf.  Section~\ref{subsec:genus:3:end}), we have $\WL(K) = \Gamma$, so the result of \cite{AGR} does not tell us the image of Weierstrass points under $\tau$. 

For small $\varepsilon > 0$, we set $B = \left[\varepsilon, \frac{1}{2} - \varepsilon\right]$. 
By \eqref{eqn:AGR:0:B}, we compute 
$
\mu_{K}(B) = -2 + 3 = 1.  
$ 
Thus $g \cdot \mu_{K}(B) = 3$. 
Since any point $p \in B$ has $\underline{n}_{K}(p) = (1, 2, 4)$, if $x \in \WL_{\mathcal{C}}(K_{\mathcal{C}})$ lands to $B$, then $\underline{n}_{K}(x) = (1, 2, 4)$ by Baker's specialization lemma. Then by 
Theorem~\ref{thm:M:m}, there are at most $3$ points $x_1, x_2, x_3 \in \WL_{\mathcal{C}}(K_{\mathcal{C}})$ that land to $(0, 1/2)$ and each of them 
has $\underline{n}_{K}(x_i) = (1, 2, 4)$. 
\end{Example}

%
\section{Proof of Theorem~\ref{thm:max:WP:genus:3}}
\label{sec:proof:thm:max:WP:genus:3}
Let $\Gamma$ be a tropical curve of genus $g \geq 2$. An edge $e \in E(\Gamma)$ is a bridge if the deletion of (the inner points of) $e$ from $\Gamma$ makes $\Gamma$ disconnected. Throughout this section, we assume that $\Gamma$ is bridgeless, i.e., $\Gamma$ has no bridges. 
Let $(G=(V,E), \ell)$ be the minimal model of $\Gamma$ obtained by 
 setting $V \colonequals \{v \in \Gamma \mid \val(v) \geq 3\}$, and defining the set of edges $E$ in bijection with the connected components of $\Gamma \setminus V$.
Let $\ell\colon E(\Gamma) \to \RR_{>0}$ denote the corresponding length function with $\ell(e)$ the length of the connected component corresponding to $e$.

Recall from the introduction that a cut vertex $v \in V(\Gamma)$ is such that the deletion of $v$ from $\Gamma$ makes $\Gamma$ disconnected. 

We assume that $g(\Gamma) = 3$ in Sections~\ref{subsec:genus:3:start}--\ref{subsec:genus:3:end}, and $g(\Gamma) = 2$ in Section~\ref{subsec:genus:2}. We give a proof of Theorem~\ref{thm:max:WP:genus:3}. 

We emphasize that within a fixed topological type, the Weierstrass locus $\WL(K)$ of a tropical curve $\Gamma$ varies in general.  For example, there exists a tropical curve $\Gamma_1$ (resp. $\Gamma_2$, resp. $\Gamma_3$) whose topological type is the complete graph $K_4$ of genus $3$ such that $\WL(K_{\Gamma_1})$ consists of $4$ points (resp. $\WL(K_{\Gamma_2})$ consists of $8$ points, resp. $\WL(K_{\Gamma_3})$ consists of $6$ points and one segment). However, Theorem~\ref{thm:max:WP:genus:3} says that 
$\sum_{\text{$p$: maximal Weierstrass point}} \wt_{K_{\Gamma_i}}(p) = 8$ for all $i = 1, 2, 3$. 

\subsection{Non-hyperelliptic curve of genus \texorpdfstring{$3$}{3} without cut vertices}
\label{subsec:genus:3:start}
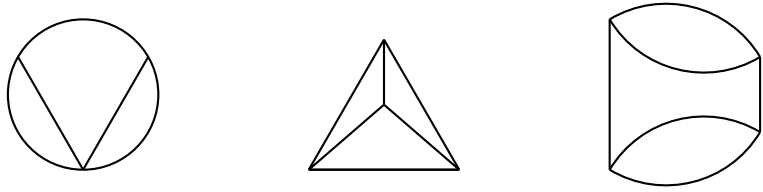
\begin{figure}[htb!] 
\[
\begin{tikzpicture}[scale=0.5]
\draw[shift={(0,2)}][thick] (0,0) circle [ radius=2];
\draw[shift={(0,2)}][thick] (0,-2) --(-1.73,1) ;
\draw[shift={(0,2)}][thick] (0,-2) --(1.73,1) ;
\draw[shift={(6,0)}] [thick] (0,0) --(4,0) ;
\draw[shift={(6,0)}] [thick] (0, 0)--(2, 3.46) ; 
\draw[shift={(6,0)}] [thick] (4, 0)--(2, 3.46) ; 
\draw[shift={(6,0)}] [thick] (0,0) --(4,0) ;
\draw[shift={(6,0)}] [thick] (0, 0)--(2, 1.73) ; 
\draw[shift={(6,0)}] [thick] (4, 0)--(2, 1.73) ; 
\draw[shift={(6,0)}] [thick] (2, 1.73)--(2, 3.46) ; 
\draw[shift={(14,0)}] [thick]  (4,3) to [out=45,in=135,relative] (0,4);
\draw[shift={(14,0)}] [thick]  (4,3) to [out=-45,in=-135,relative] (0,4);
\draw[shift={(14,0)}] [thick]  (4,1) to [out=45,in=135,relative] (0,0);
\draw[shift={(14,0)}] [thick]  (4,1) to [out=-45,in=-135,relative] (0,0);
\draw[shift={(14,0)}] [thick] (0,0) --(0,4) ;
\draw[shift={(14,0)}] [thick] (4,1) --(4,3) ;
\end{tikzpicture}
\]
\caption{Non-hyperelliptic curves of genus $3$ without cut vertices.}
\label{fig:subsec:genus:3:start}
\end{figure}

Let $\Gamma$ be a non-hyperelliptic curve of genus $3$ without cut vertices. Then $\Gamma$ is one of the three tropical curves in Figure~\ref{fig:subsec:genus:3:start}. 
(Here, for the rightmost tropical curve, the lengths of the two vertical edges are not equal to each other.) 
We are going to show that $\sum_{\text{$p$: maximal Weierstrass point}} \wt_{K_{\Gamma_i}}(p) = 8$. 

Let $p \in \Gamma$ be a Weierstrass point, i.e., $\rank(3 (p)) \geq 1$. 
By the tropical Riemann--Roch theorem (see Theorem~\ref{thm:RR}), we have 
$\rank(K - 3(p)) = \rank(3 (p)) - 1 \geq 0$, 
so there exists $q \in \Gamma$ with $K \sim 3(p) + (q)$. Since 
$\Gamma$ has no bridges, $q$ is uniquely determined by $p$. 

Note that there is no $p \in \Gamma$ 
with $\underline{n}_{K}(p) = (1, 3, 5)$, because 
$\Gamma$ is non-hyperelliptic. 
For $p \in \Gamma$, we see that 
$\underline{n}_{K}(p) = (1, 2, 3)$ if and only if $\rank(K - 3(p)) = -1$, and 
\begin{align*}
\text{$\underline{n}_{K}(p) = (1, 2, 4)$}  
& \quad\Longleftrightarrow\quad
\text{there exists a unique $q \in \Gamma$ with $K \sim 3(p) + (q)$ with $q \neq p$},   
\\
\text{$\underline{n}_{K}(p) = (1, 2, 5)$} 
& \quad\Longleftrightarrow\quad
\text{$K \sim 4(p)$}.  
\end{align*}

Let $A$ be a connected component of the Weierstrass locus of $\Gamma$. 
We assume that $A$ is not a singleton. 

\begin{Lemma}
\label{lem:nonhyp:3:1}
There is an edge $e \in E(\Gamma)$ such that $A \subseteq e$.   
\end{Lemma}

\Proof
Let $p_0 \in \partial A$.  Since $A$ is not a singleton, there is a tangent direction $\nu \in T_{p_0}\Gamma$ and an edge $e \in E(\Gamma)$ with $p_0 \in e$ such that $e$ is identified with the closed interval $[0, \ell(e)]$ of length $\ell(e)$ with the direction $0$ to $\ell(e)$ agreeing with $\nu$ and that $p_0 + \varepsilon \in A$ for any small positive $\varepsilon$, where we identify $p_0$ with its coordinate in $[0, \ell(e)]$.  

We fix a small positive $\varepsilon$, and write $p_x = p_0 + x$ for $0 \leq x \leq \varepsilon$. Since $p_x \in A$, there exists a unique $q_x \in \Gamma$ 
with $K \sim 3(p_x) + (q_x)$. 

\begin{Claim}
$\underline{n}_{K}(p_0) = (1, 2, 4)$. 
\end{Claim}

Indeed, since $A$ is closed (see Theorem~\ref{thm:upper}) and $p_0 \in A$, we have $\underline{n}_{K}(p_0) = (1, 2, 4)$ or $(1, 2, 5)$. 
To derive a contradiction, assume that $\underline{n}_{K}(p_0) = (1, 2, 5)$. Then $K \sim 4 (p_0)$. Suppose that $p_{0} \in V(\Gamma)$. Then $\val(p_0) \geq 3$. We denote the edges with endpoints $p_0$ by $e, e_1, \ldots, e_N$ (with $N \colonequals \val(p_0) - 1 \geq 2$). Since $\Gamma \setminus \{p_0\}$ is connected by our assumption that $\Gamma$ has no cut vertices, for small $x > 0$, for each $1 \leq i \leq N$ there is a path $\gamma_i$ that connects $p_x$ ($x > 0$) and a point on $e_i$ without passing through $p_0$. Since $N \geq 2$, it follows that $K \not\sim 3 (p_x) + (q_x)$ for small $x > 0$ by e.g. Luo's algorithm \cite{Luo}, a contradiction. Next, suppose that $p_{0} \not\in V(\Gamma)$. Then $p_0$ is an inner point of $e$. 
We see that there is a small positive $\varepsilon^\prime$ such that for any $- \varepsilon^\prime \leq x \leq \varepsilon^\prime$, we have 
$K \sim 4 (p_0) \sim 3(p_0 + x) + (p_0 - 3x)$. It follows that there is 
an neighborhood $U_0$ of $p_0$ such that every point $p \in U_0$ is a Weierstrass point. 
This contradicts with $p_0 \in \partial A$. Thus we get $\underline{n}_{K}(p_0) = (1, 2, 4)$, and we have shown the claim. 

We replace $\varepsilon$ by a smaller number if necessary, we may assume that $\val(p_x) = 2$ and $\val(q_x) = 2$ for any $0 < x \leq \varepsilon$. 
For any $0 < x_1 < x_2 \leq \varepsilon$, since $3(p_x) + (q_x)$ is $(p_x)$-reduced, we see that $p_x$ moves continuously from $p_{x_1}$ to $p_{x_2}$ and $q_x$ moves continuously from $q_{x_1}$ to $q_{x_2}$ by, for example, Luo's algorithm \cite{Luo}. 
Shrinking $\varepsilon$ if necessary, we may assume that the moving path  from $p_{x_1}$ to $p_{x_2}$ and that from $q_{x_1}$ to $q_{x_2}$ are disjoint. 
We identify the moving path from $p_{x_1}$ to $p_{x_2}$ with the interval $[x_1, x_2]$ and that from $q_{x_1}$ to $q_{x_2}$ with the interval $[y_1, y_2]$. 

It follows that there exists $f_{12} \in \Rat(\Gamma)$ such that
$3(p_{x_1}) + (q_{x_1})  + \zero(f_{12}) = 3(p_{x_2}) + (q_{x_2})$ 
with $f_{12}(x) = 3 (x-x_1)$ for $x \in [x_1, x_2]$, $f_{12}(y) = (y-y_1)$ and $f_{12}$ is locally constant on $\Gamma\setminus [x_1, x_2] \cup [y_1, y_2]$. 
In the right tropical curve in the above figure, since $p_{x_1}$ is a Weierstrass point, we see that $p_{x_1}$ does not lie on the inner points of the vertical edges. It follows from the above figures, if $q_{y_1}$ is not on $e$, then $\Gamma\setminus [x_1, x_2] \cup [y_1, y_2]$ is connected, and there does not exists such an $f_{12}$, because $f_{12}(x_1) \neq f_{12}(x_2)$. 

We conclude that $q_{y_1} \in e$. Moving $x_1$ to $0$, we see that $p_0, q_0 \in e$. Since $p_x \in A$, we get $0 \leq p_0 < q_0 \leq \ell(e)$. 

\medskip

{\bf Case 1}: Suppose that $p_0 = 0$. Then $0 < q_0 \leq \ell(e)$. 
We have $A =\left[0, \frac{q_0}{3}\right] \subseteq e$.  

{\bf Case 2}: Suppose that $0 < p_0 < \frac{2 \ell(e)}{3}$. Since $p_0 \in \partial A$, we have  $q_0 = \ell(e)$. 
We get $A =\left[p_0, p_0 +\frac{\ell(e)}{3}\right] \subseteq e$. 

{\bf Case 3}: Suppose that $\frac{2 \ell(e)}{3} \leq p_0 < \ell(e)$. Since $p_0 \in \partial A$, we have  $q_0 = \ell(e)$. 
We get $A =\left[p_0, \ell(e)\right] \subseteq e$. 

In all cases, we have $A \subseteq e$. 
\QED

\begin{Lemma}
\label{lem:nonhyp:3:2}
On $A$, there exists a unique point $p$ with $\underline{n}_{K}(p) = (1, 2, 5)$. 
\end{Lemma}

\Proof
In Case 1 in the proof of Lemma~\ref{lem:nonhyp:3:1}, 
$p = \frac{1}{4} q_0$ is a unique point on $A$ with $\underline{n}_{K}(p) = (1, 2, 5)$. Similarly, in Cases~2 and~3, $p = \frac{3}{4} p_0 + \frac{1}{4} \ell(e)$ is a unique point on $A$ with $\underline{n}_{K}(p) = (1, 2, 5)$. 
\QED

\begin{Lemma}
\label{lem:nonhyp:3:3}
$\mu_{K}(A) = 2$. 
\end{Lemma}

\Proof
Consider Case 1 in the proof of Lemma~\ref{lem:nonhyp:3:1}. 
We take $f_0 \in \Rat(\Gamma)$ with $K + \zero(f_0) = 3(p_0)+ (q_0) $, where $p_0 = 0$. We denote the edges with endpoints $p_0$ by $e, e_1, \ldots, e_N$ (with $N \colonequals \val(p_0)-1\geq 2$). Let $\nu, \nu_1, \ldots, \nu_N$ be tangent directions emanating from $p_0 = 0$ on $e, e_1, \ldots, e_N$. 
We write $\slope_{p_0, \nu}(f_0)= s_0$ and $\slope_{p_0, \nu_i}(f_0)= s_0^{(i)}$. Since the coefficient of $K$ at $p_0$ is $(N+1) - 2 = N-1$, we see that $s_0 + s_0^{(1)} + \cdots + s_0^{(N)} = 
-\zero(f_0)(p_0) = - 3 + (N-1) = N-4$. 

Set 
\[
g(x)
= 
\begin{cases}
2x &  x \in \left[0, \frac{1}{3} q_0\right],  \\
- x + q_0 & x \in [\frac{1}{3} q_0, q_0], \\
0 & \text{elsewhere}. 
\end{cases}
\]
Then $K + \zero(f_0) + \zero(g) = 3 \left(\frac{1}{3} q_0\right) + (0)$. 

Let $\nu$ (the same notation) be the tangent direction at $\frac{1}{3} q_0$ from $0$ to $\ell(e)$. We get 
 $\slope_{\nu, \frac{1}{3} q_0}(f_0 + g)= s_0 - 1$. 
 
It follows from \eqref{eqn:AGR} that 
\begin{align*}
\mu_{K}(A) 
& = (\val(p_0) - 2) 
+ 
(g(A)-1)(g-1) - (s_0^{(1)} + \cdots + s_0^{(N)}) -  \slope_{\nu, \frac{1}{3} q_0}(f_0 + g) \\
& = (N-1) -2 - (s_0^{(1)} + \cdots + s_0^{(N)}) - (s_0 - 1) = 2. 
\end{align*}

\medskip
Consider Case 2. 
We take $f_0 \in \Rat(\Gamma)$ with $K + \zero(f_0) = 3(p_0)+ (q_0) $, where $q_0 = \ell(e)$. Let $\nu$ be the tangent direction from $0$ to $\ell(e)$ on $e$, and 
$\nu^\prime$ be the opposite direction of $\nu$. 
We write $\slope_{p_0, \nu}(f_0)= s_0$ and $\slope_{p_0, \nu^\prime}(f_0)= s_0^\prime$. Then $s_0 + s_0^\prime = -3$. Since the support of $K$ is included in $V(\Gamma)$, $f_0$ has constant slope on $[p_0, \ell(e)]$. 

Set 
\[
g(x)
= 
\begin{cases}
-x + p_0 & x \in [0, p_0], \\
2(x-p_0) &  x \in \left[p_0, p_0 + \frac{\ell(e)}{3}\right],  \\
- x +p_0  + \ell(e) & x \in [p_0 + \frac{\ell(e)}{3}, \ell(e)], \\
p_0 & x \in \Gamma\setminus e. 
\end{cases}
\]
Then $K + \zero(f_0) + \zero(g) = 3 \left(p_0 + \frac{\ell(e)}{3}\right) + (0)$. 

Let $\nu$ (the same notation) be the tangent direction at $p_0 + \frac{\ell(e)}{3}$ from $0$ to $\ell(e)$, and $\nu^\prime$ the opposite direction. We get 
 $\slope_{\nu, p_0+ \frac{\ell(e)}{3}}(f_0 + g)= s_0 - 1$. 

It follows that 
\[
\mu(A) = (g(A)-1)(g-1) - \slope_{p_0, \nu^\prime}(f_0) -  \slope_{p_0 + \frac{\ell(e)}{3}, \nu}(f_0+g) 
= -2 - s_0^\prime - (s_0 - 1) = 2. 
\]

\medskip
Case 3 is similar to Case 1. 
\QED

\bigskip
{\sl Proof of Theorem~\ref{thm:max:WP:genus:3}~(1)}: \quad
Let $A_1, \ldots, A_m$ be the set of Weierstrass locus of $\Gamma$ such that 
$A_1, \ldots, A_k$ are singletons, and $A_{k+1}, \ldots, A_m$ are not singletons for some $0 \leq k \leq m$.  We put $p_i = A_i$ for $1 \leq i \leq k$. 
For $k+1 \leq i \leq m$, by Lemma~\ref{lem:nonhyp:3:2}, 
let $p_i$ be the unique point on $A_i$ with 
$\underline{n}(p_i) = (1, 2, 5)$. Then $\{p_1, \ldots p_k, p_{k+1}, \ldots, p_m\}$
is the set of maximal Weierstrass points of $\Gamma$. 
For $k+1 \leq i \leq m$, Lemma~\ref{lem:nonhyp:3:3} tells us that  
$\mu(A_i) = 2 = \wt_\Gamma(p_i)$. Then Theorem~\ref{thm:AGR} implies that \[
\sum_{p_i\colon \text{maximal Weierstrass point}}
\wt_{K}(p_i)
= \sum_{i=1}^m \mu(A_i) = g^2 - 1 = 8. 
\]
Since $\mu(A_i)$ is either $1$ or $2$, we see that $4 \leq m \leq 8$. 
\QED

\subsection{Hyperelliptic curve of genus \texorpdfstring{$3$}{3} without cut vertices}
\label{subsec:hyperell:wituout:cuts}
\begin{figure}[htb!]
\[
\begin{tikzpicture}[scale=0.5]
\draw[shift={(0,1)}][thick] (4,0) to [out=20,in=160,relative] (0,0);
\draw[shift={(0,1)}][thick] (4,0) to [out=-20,in=-160,relative] (0,0);
\draw[shift={(0,1)}][thick] (4,0) to [out=80,in=100,relative] (0,0);
\draw[shift={(0,1)}][thick] (4,0) to [out=-80,in=-100,relative] (0,0);
\draw[shift={(10,0)}] [thick]  (4,2) to [out=45,in=135,relative] (0,2);
\draw[shift={(10,0)}] [thick]  (4,2) to [out=-45,in=-135,relative] (0,2);
\draw[shift={(10,0)}] [thick]  (4,0) to [out=45,in=135,relative] (0,0);
\draw[shift={(10,0)}] [thick]  (4,0) to [out=-45,in=-135,relative] (0,0);
\draw[shift={(10,0)}] [thick] (0,0) --(0,2) ;
\draw[shift={(10,0)}] [thick] (4,0) --(4,2) ;
\end{tikzpicture}
\]
\caption{Hyperelliptic curve of genus $3$ without cut vertices.}
\label{fig:hyperell:wituout:cuts}
\end{figure}
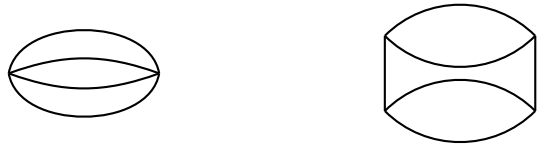

Let $\Gamma$ be a hyperelliptic curve of genus $3$ without cut vertices. 
Then $\Gamma$ is one of the tropical curves in Figure~\ref{fig:hyperell:wituout:cuts}. 

(Here, for the right tropical curve, the lengths of the two vertical edges are  equal to each other.)

\medskip
{\sl Proof of Theorem~\ref{thm:max:WP:genus:3}~(2)}: \quad
For the left tropical curve $\Gamma$, let $\{e_1, \ldots, e_4\}$ be the set of edges. Let $p_i$ be the middle point of each $e_i$. The computation in Example~\ref{eg:dipole:3} show that $\{p_1, \ldots, p_4\}$ is 
the set of maximal Weierstrass points of $\Gamma$ and $\wt_\Gamma(p_i) = \frac{g(g-1)}{2} = 3$. It follows that 
\[
\sum_{p_i\colon \text{maximal Weierstrass point}}
\wt_{K}(p_i)
= (g+1) \frac{g(g-1)}{2} = 12.  
\] 

For the right tropical curve $\Gamma$, let $\{e_1, \ldots, e_4\}$ be the set of horizontal edges, and let $p_i$ be the middle point of each $e_i$. Then similar to the computation in Example~\ref{eg:dipole:3}, we see that $\{p_1, \ldots, p_4\}$ is the set of maximal Weierstrass points of $\Gamma$, and 
the sum of the total weight of these points is equal to $12$. 
\QED

\subsection{Non-hyperelliptic curve of genus \texorpdfstring{$3$}{3} with cut vertices}
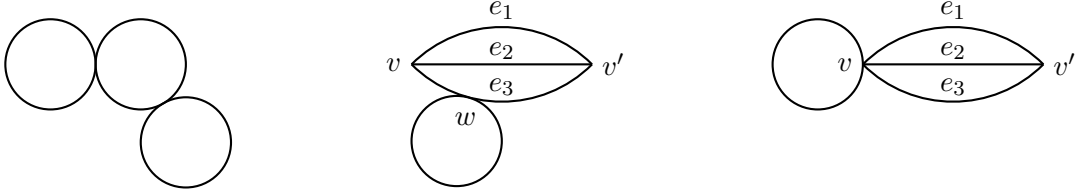
\begin{figure}[htb!]
\[
\begin{tikzpicture}[scale=0.6]
\draw[shift={(0,0)}] [thick] (0,0) circle [ radius=1];
\draw[shift={(0,0)}] [thick] (2,0) circle [ radius=1 ];
\draw[shift={(0,0)}] [thick] (3,-1.73) circle [ radius=1 ];
\draw[shift={(8,0)}] [thick] (0,0) --(4,0) ;
\draw[shift={(8,0)}] [thick]  (4,0) to [out=45,in=135,relative] (0,0);
\draw[shift={(8,0)}] [thick]  (4,0) to [out=-45,in=-135,relative] (0,0);
\draw[shift={(8,0)}] [thick] (1,-1.7) circle [ radius=1 ];
\node at (10,1.2) {$e_1$};
\node at (10,0.3) {$e_2$};
\node at (10,-0.5) {$e_3$};
\node[left] at (8, 0) {$v$};
\node[right] at (12, 0) {$v^\prime$};
\node at (9.2, -1.2) {$w$};
\draw[shift={(16,0)}] [thick] (2,0) --(6,0) ;
\draw[shift={(16,0)}] [thick]  (6,0) to [out=45,in=135,relative] (2,0);
\draw[shift={(16,0)}] [thick]  (6,0) to [out=-45,in=-135,relative] (2,0);
\draw[shift={(16,0)}] [thick] (1,0) circle [ radius=1 ];
\node at (20,1.2) {$e_1$};
\node at (20,0.3) {$e_2$};
\node at (20,-0.5) {$e_3$};
\node[left] at (18, 0) {$v$};
\node[right] at (22, 0) {$v^\prime$};
\end{tikzpicture}
\]
\caption{Non-hyperelliptic curve of genus $3$ with cut vertices}
\label{fig:non:hyperell:with:cuts}
\end{figure}

Let $\Gamma$ be a non-hyperelliptic curve of genus $3$ with cut vertices. 
Then $\Gamma$ is one of the tropical curves in Figure~\ref{fig:non:hyperell:with:cuts}. 

(Here, for the left tropical curve, the two vertices are not antipodal with respect to the middle circle. For the middle tropical curve, the vertex $w$ on the circle is not the middle point of the edge $e_3$.) For these graphs, it is easy to compute the Weierstrass gap sequences. 

We write the left, middle, and right tropical curves for $\Gamma_1, \Gamma_2$, and $\Gamma_3$. Black, thick gray, and blue points are respectively points with Weierstrass gap sequence $(1, 2, 3)$, $(1, 2, 4)$ and $(1, 2, 5)$. 

\begin{itemize}
\item
For $\Gamma_1$, the computation of Weierstrass gap sequences is easy, 
and we have the following figure. 
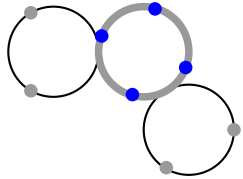
\[
\begin{tikzpicture}[scale=0.6]
\draw [thick] (0,0) circle [ radius=1];
\draw [thick] (3,-1.73) circle [ radius=1 ];
\draw[color=black!40, line width=1mm] (2,0) circle [ radius=1 ];
\foreach \c in {
(-0.5,1.73/2), (-0.5,-1.73/2), (2.5, -2.57), (4, -1.73)
} {
	\fill[color=black!40] \c circle [radius=0.15];
}
\foreach \c in {
(1.75, -0.95), (2.25, 0.95), (1.07, 0.35), (2.93, -0.35)
} {
	\fill[color=blue] \c circle [radius=0.15];
}
\end{tikzpicture}
\]
We get 
\[
\sum_{\text{$p$: maximal Weierstrass point of $\Gamma_1$}}
 \wt_{\Gamma}(p) = 8+4 = 12. 
\]
\item
For $\Gamma_2$, we argue as in Section~\ref{subsec:genus:3:start}. 
Let $v, v^\prime, w$ be the vertices as in the above figure, 
and let $e_1, e_2, e_3$ be the edges connecting $v$ and $v^\prime$. 
We identify $e_3$ with the closed interval $[0, \ell(e)]$ with $v, v^\prime$ corresponding to $0, \ell(e_3)$, respectively. Let $a$ ($0 < a < \ell(e_3)$) be the coordinate of $w$. We may and do assume that $a < \ell(e_3)/2$. 

Let $A$ be a connected component of the Weierstrass locus of $\Gamma_2$. 
Since $K = 2(w) + (v) + (v^\prime)$, if $w \in A$, we see that 
$A = \left[\frac{2}{3}a, a + \frac{1}{3} \ell(e_3)\right]$. We also have 
$\mu(A) = 2$ and there exists a unique point $q = \frac{3}{4}a + \frac{1}{4} \ell(e_3) \in A$ with $\underline{n}_{\Gamma_2}(q) = (1, 2, 5)$. 

If $A$ is not a singleton and $w \not\in A$, then the proofs of Lemmas~\ref{lem:nonhyp:3:1} to~\ref{lem:nonhyp:3:3} show that 
$A \subseteq e_1$ or $A \subseteq e_2$, and that $\mu(A) = 2$ and 
$q \in A$ with $\underline{n}_{\Gamma_2}(q) = (1, 2, 5)$. 

By Theorem~\ref{thm:AGR} and arguing as in Section~\ref{subsec:genus:3:start}, there are in total $m$ maximal Weierstrass points on $\Gamma_2$ for some $4 \leq m \leq 8$, and 
we have 
\[
\sum_{\text{$p$: maximal Weierstrass point of $\Gamma_2$}}
 \wt_{\Gamma}(p) = 8. 
\]

\item
For $\Gamma_3$, we can argue as for $\Gamma_2$, but here we will compute Weierstrass gap sequences directly. 
Let $v, v^\prime$ be the vertices as in the above figure, 
and let $e_1, e_2, e_3$ be the edges connecting $v$ and $v^\prime$. 
We identify $e_i$ with the closed interval $[0, \ell(e_i)]$ with $v, v^\prime$ corresponding to $0, \ell(e_i)$, respectively. Without loss of generality, 
we may assume that $\ell(e_1) \leq \ell(e_2) \leq \ell(e_3)$. Then we get the following. 
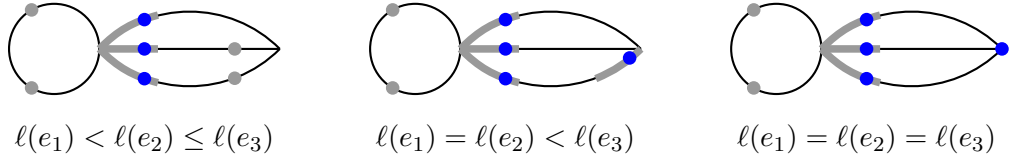
\[
\begin{tikzpicture}[scale=0.6]
\begin{scope}
\draw [thick] (2,0) --(6,0) ;
\draw [thick]  (6,0) to [out=45,in=135,relative] (2,0);
\draw [thick]  (6,0) to [out=-45,in=-135,relative] (2,0);
\draw [thick] (1,0) circle [ radius=1 ];
\foreach \c in {
(1/2, 1.73/2), (1/2, -1.73/2), (5, 0), (5, -0.65)
} {
	\fill[color=black!40] \c circle [radius=0.15];
}
\draw[color=black!40] [line width = 1mm]  (2, 0) to [out=15,in=165,relative] (3.3, 0.75);
\draw[color=black!40] [line width = 1mm]  (2, 0) to [out=-15,in=-165,relative] (3.3, -0.75);
\draw[color=black!40] [line width = 1mm]  (2, 0)--(3.3, 0);
\foreach \c in {
(3, 0.65), (3, 0), (3, -0.66)
} {
	\fill[color=blue] \c circle [radius=0.15];
}
\end{scope}
\node at (3,-2) {$\ell(e_1) < \ell(e_2) \leq \ell(e_3)$};
\begin{scope}[shift={(8,0)}]
\draw [thick] (2,0) --(6,0) ;
\draw [thick]  (6,0) to [out=45,in=135,relative] (2,0);
\draw [thick]  (6,0) to [out=-45,in=-135,relative] (2,0);
\draw [thick] (1,0) circle [ radius=1 ];
\draw[color=black!40] [line width=1mm]  (2, 0) to [out=15,in=165,relative] (3.3, 0.75);
\draw[color=black!40] [line width=1mm]  (2, 0) to [out=-15,in=-165,relative] (3.3, -0.75);
\draw[color=black!40] [line width=1mm]  (2, 0)--(3.3, 0);
\foreach \c in {
(1/2, 1.73/2), (1/2, -1.73/2)
} {
	\fill[color=black!40] \c circle [radius=0.15];
}
\draw[color=black!40] [line width=1mm]  (5, -0.66) to [out=-10,in=-170,relative] (6, 0);
\foreach \c in {
(3, 0.65), (3, 0), (3, -0.66), (5.75, -0.2)
} {
	\fill[color=blue] \c circle [radius=0.15];
}
\end{scope}
\node at (11,-2) {$\ell(e_1) = \ell(e_2) < \ell(e_3)$};
\begin{scope}[shift={(16,0)}]
\draw [thick] (2,0) --(6,0) ;
\draw [thick]  (6,0) to [out=45,in=135,relative] (2,0);
\draw [thick]  (6,0) to [out=-45,in=-135,relative] (2,0);
\draw [thick] (1,0) circle [ radius=1 ];
\draw[color=black!40, line width=1mm]  (2, 0) to [out=15,in=165,relative] (3.3, 0.75);
\draw[color=black!40, line width=1mm]  (2, 0) to [out=-15,in=-165,relative] (3.3, -0.75);
\draw[color=black!40, line width=1mm]  (2, 0)--(3.3, 0);
\foreach \c in {
(1/2, 1.73/2), (1/2, -1.73/2)
} {
	\fill[shift={(0,0)}, color=black!40] \c circle [radius=0.15];
}
\foreach \c in {
        (3, 0.65), (3, 0), (3, -0.66), (6, 0)
    } {
        \fill[color=blue] \c circle [radius=0.15];
    }
\end{scope}
\node at (19,-2) {$\ell(e_1) = \ell(e_2) = \ell(e_3)$};
\end{tikzpicture}
\]
(In the middle figure, the coordinates of the blue points on $e_3$ are 
$\frac{1}{4} \ell(e_3)$ and $\frac{3}{4} \ell(e_3) + \frac{1}{4} \ell(e_1)$. With some lengths of $e_i$, it is possible that every point on $e_3$ is a blue point.)
In any case, we have 
\[
\sum_{\text{$p$: maximal Weierstrass point of $\Gamma_3$}}
 \wt_{\Gamma}(p) = 10. 
\]
\end{itemize}

\subsection{Hyperelliptic curve of genus \texorpdfstring{$3$}{3} with cut vertices}
\label{subsec:genus:3:end}
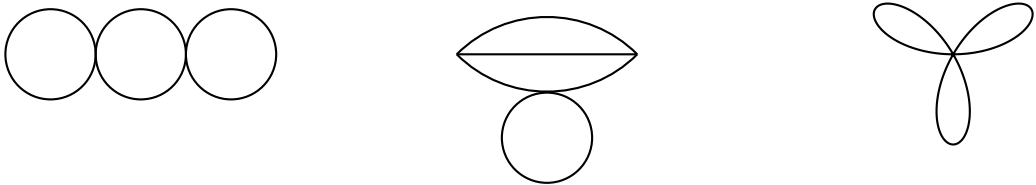
\begin{figure}[htb!]
\[
\begin{tikzpicture}[scale=0.6]
\draw[shift={(0,0)}] [thick] (0,0) circle [ radius=1];
\draw[shift={(0,0)}] [thick] (2,0) circle [ radius=1 ];
\draw[shift={(0,0)}] [thick] (4,0) circle [ radius=1 ];
\draw[shift={(9,0)}] [thick] (0,0) --(4,0) ;
\draw[shift={(9,0)}] [thick]  (4,0) to [out=45,in=135,relative] (0,0);
\draw[shift={(9,0)}] [thick]  (4,0) to [out=-45,in=-135,relative] (0,0);
\draw[shift={(9,0)}] [thick] (2,-1.85) circle [ radius=1 ];
\draw[shift={(20,0)}][thick, black, domain=0:180, samples=200] 
        plot (\x:{2*sin(3*\x)});
\end{tikzpicture}
\]
\caption{Hyperelliptic curve of genus $3$ with cut vertices}
\label{fig:genus:3:end}
\end{figure}

Let $\Gamma$ be a hyperelliptic curve of genus $3$ with cut vertices. 
Then $\Gamma$ is one of the tropical curves in Figure~\ref{fig:genus:3:end}. 

(Here, for the left tropical curve, the two vertices are antipodal with respect to the middle circle. For the right tropical curve, the vertex on the circle is the middle point of the edge.) For these graphs, it is easy to compute the Weierstrass gap sequences. 

Black, thick gray, blue, and green-star points are respectively points with Weierstrass gap sequence $(1, 2, 3)$, $(1, 2, 4)$, $(1, 2, 5)$ and $(1, 3, 5)$. Then we have the following. 
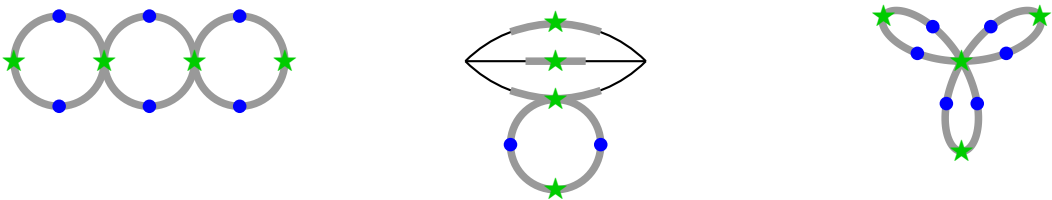
\[
\begin{tikzpicture}[scale=0.6]
\draw[shift={(0,0)}, color=black!40] [line width=1mm] (0,0) circle [ radius=1];
\draw[shift={(0,0)}, color=black!40] [line width=1mm] (2,0) circle [ radius=1 ];
\draw[shift={(0,0)}, color=black!40] [line width=1mm] (4,0) circle [ radius=1 ];
\foreach \c in {
(-1, 0), (1, 0), (3, 0), (5, 0)
} {
	\node[green!80!black, font=\small] at \c {$\bigstar$};
}
\foreach \c in {
(0, 1), (0, -1), (2, 1), (2, -1), (4, 1), (4, -1)
} {
	\fill[color=blue] \c circle [radius=0.15];
}
\begin{scope}[shift={(9,0)}]
\draw [thick] (0,0) --(4,0) ;
\draw [thick]  (4,0) to [out=45,in=135,relative] (0,0);
\draw [thick]  (4,0) to [out=-45,in=-135,relative] (0,0);
\draw[color=black!40] [line width=1mm] (2,-1.85) circle [ radius=1 ];
\draw[color=black!40] [line width=1mm] (4/3,0) --(8/3,0) ;
\draw[color=black!40] [line width=1mm]  (1, 0.65) to [out=18,in=162,relative] (3, 0.65);
\draw[color=black!40] [line width=1mm]  (1, -0.65) to [out=-18,in=-162,relative] (3, -0.65);
\foreach \c in {
(2, 0), (2, 0.85), (2, -0.85), (2, -2.85)
} {
	\node[green!80!black, font=\small] at \c {$\bigstar$};
}
\foreach \c in {
(1, -1.85), (3, -1.85)
} {
	\fill[color=blue] \c circle [radius=0.15];
}
\end{scope}
\begin{scope}[shift={(20,0)}]
\draw[line width=1mm, black!40, domain=0:180, samples=200] 
        plot (\x:{2*sin(3*\x)});
\node[green!80!black, font=\small] at (0,0) {$\bigstar$};
\foreach \angle in {30, 90, 150} {
            \node[green!80!black, font=\small] at (\angle:{2*sin(3*\angle)}) {$\bigstar$};
        }
\foreach \angle in {10, 50, 70, 110, 130, 170} {
            \fill[blue] (\angle:{2*sin(3*\angle)}) circle [radius=0.15];
        }
\end{scope}
\end{tikzpicture}
\]

We write the left, middle, and right tropical curves for $\Gamma_4, \Gamma_5$, and 
$\Gamma_6$. Then we get 
\begin{align*}
\sum_{p}
 \wt_{\Gamma_4}(p) = 24, 
 \qquad
\sum_{p}
 \wt_{\Gamma_5}(p) = 16, 
 \qquad
 \sum_{p}
 \wt_{\Gamma_6}(p) = 24, 
\end{align*}
where the $p$ run through all the maximal Weierstrass points of $\Gamma_i$ ($i = 4, 5, 6$). 

\subsection{Genus \texorpdfstring{$2$}{2} tropical curves}
\label{subsec:genus:2}
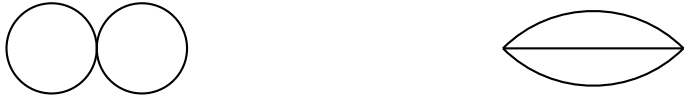
\begin{figure}[htb!]
\[
\begin{tikzpicture}[scale=0.6]
\draw[shift={(0,0)}] [thick] (0,0) circle [ radius=1];
\draw[shift={(0,0)}] [thick] (2,0) circle [ radius=1 ];
\draw[shift={(10,0)}] [thick] (0,0) --(4,0) ;
\draw[shift={(10,0)}] [thick]  (4,0) to [out=45,in=135,relative] (0,0);
\draw[shift={(10,0)}] [thick]  (4,0) to [out=-45,in=-135,relative] (0,0);
\end{tikzpicture}
\]
\caption{Genus $2$ tropical curves}
\label{fig:genus:2}
\end{figure}

Let $\Gamma$ be a tropical curve of genus $2$ without bridges. 
Then $\Gamma$ is one of the tropical curves in Figure~\ref{fig:genus:2}. 

The right graph is a dipole graph $B_2$ of genus $2$ (see Example~\ref{eg:dipole}). Let $q_1, q_2, q_3$ be the middle points of three edges of $B_2$. Then we see that $\underline{n}_{K}(p) = (1, 2)$ if $p \neq q_i$ and $\underline{n}_{K}(q_i) = (1, 3)$ for $1 \leq i \leq 3$. For the left graph, 
we compute similarly.  
We get 
\[
\sum_{\text{$p$: maximal Weierstrass point}}
 \wt_{\Gamma}(p) = 6.
\]

For a related question, see Question~\ref{q:max:W:1}.

%
\section{Proof of Theorem~\ref{thm:main:1}}
\label{sec:classification}

Recall from the introduction (see \eqref{eqn:def:S:g:canonical}) that 
\[
\mathcal{S}(g) = \{\underline{n} =  (n_1, \ldots, n_g) \in \ZZ^g \mid 1 \leq n_1 < \cdots < n_g \leq 2g-1\}. 
\]
Recall also that $\mathcal{G}(g)$ is the subset of $\mathcal{S}(g)$ consisting of those $\underline{n}$ that are achieved by some point of some tropical curve of genus $g$, and that $\mathcal{N}(g)$ is the subset of $\mathcal{S}(g)$ containing those with the complement $\underline{n}^c$ in $\ZZ_{> 0}$ forming a numerical semigroup. 

To exclude some sequences in $\mathcal{S}(g)$ that do not belong 
to $\mathcal{G}(g)$, we use the following lemma. 

\begin{Lemma}
\label{lemma:useful:2}
Let $\Gamma$ be a tropical curve of genus $g \geq 2$, and let $p \in \Gamma$. 
Suppose that $\underline{n}_{K}(p) = (1, n_2, n_3, \ldots n_{g})$ with $3 \leq n_2$. 
Then $\Gamma$ is hyperelliptic, and $\underline{n}_{K}(p) = (1, 3, 5, \ldots, 2g-1)$. 
\end{Lemma}

\Proof
If $\underline{n}_{K}(p) = (1, n_2, n_3, \ldots n_{g})$ with $3 \leq n_2$, then $\rank(2(p)) = 1$. It follows that $\Gamma$ is hyperelliptic, and 
we have $\rank(k(p)) = \lfloor k/2\rfloor$ for $k = 0, 1, \ldots, 2g-1$. 
Then $\underline{n}_{K}(p) = (1, 3, 5, \ldots, 2g-1)$. 
\QED

When $\Gamma$ is hyperelliptic, we let $v_0 \in \Gamma$ denote a point invariant under a hyperelliptic involution (after contracting leaf edges) as in Section~\ref{subsec:trop:hyper}. We also use the following lemma. 

\begin{Lemma}
\label{lem:cor:Clifford}
Let $\Gamma$ be a hyperelliptic graph of genus $g \geq 4$, and let $D \in \Div(\Gamma)$ be a divisor 
such that $2 \leq \rank(D) \leq g-2$ and $\rank(D) = \frac{1}{2} \deg(D)$. Then 
$D \sim 2 \rank(D) (v_0)$. 
\end{Lemma}

\Proof
Since $\rank(D) \geq 2$, we may assume that $D$ is effective. 
We use Proposition~\ref{prop:rank:hyp:ell:KY}.
With the notation there, 
we claim that $\deg(D) - p_{\Gamma}(D) \leq g$. Indeed, if $\deg(D) - p_{\Gamma}(D) \geq g+1$, then, 
by Proposition~\ref{prop:rank:hyp:ell:KY}, 
$\rank(D) = \deg(D) -g = 2 \rank(D) - g$. It follows that $\rank(D) = g$, 
which contradicts with $\rank(D) \leq g-2$, and the claim follows.

Then, by Proposition~\ref{prop:rank:hyp:ell:KY}, $\rank(D) = p_{\Gamma}(D)$. By definition of $p_{\Gamma}(D)$, we have $\left|D - 2 \rank(D) (v_0)\right| \neq \emptyset$. Since $\deg(D) = 2 \rank(D)$ by assumption, we have $D \sim 2 \rank(D) (v_0)$. 
\QED

We begin the proof of Theorem~\ref{thm:main:1}. 
\subsection*{Genus \texorpdfstring{$2$}{2} case}
We have 
\[
\mathcal{S}(2) = \mathcal{N}(2) = \left\{
(1, 2), (1, 3)
\right\}. 
\]
Let $B_2$ be the dipole graph of genus $2$, where all the edges are assigned length one (see Example~\ref{eg:dipole}). By Section~\ref{subsec:genus:2}, there exist $p, q \in B_2$ with 
$\underline{n}_{K}(p) = (1, 2)$ and $\underline{n}_{K}(q) = (1, 3)$. 
We obtain Theorem~\ref{thm:main:1}~(1). 

\subsection*{Genus \texorpdfstring{$3$}{3} case}
We have 
\begin{align*}
\mathcal{S}(3) 
& = \left\{
(1, 2, 3), (1, 2, 4), (1, 2, 5), (1, 3, 4), (1, 3, 5), (1, 4, 5) 
\right\}, \\
\mathcal{N}(3) 
& = \left\{
(1, 2, 3), (1, 2, 4), (1, 2, 5), (1, 3, 5)
\right\}. 
\end{align*}

By Lemma~\ref{lemma:useful:2}, 
neither $(1, 3, 4)$ nor $(1, 4, 5)$ is of the form $\underline{n}_{K}(p)$ for a point $p$ on a tropical curve $\Gamma$ of genus $3$. Thus $\mathcal{G}(3)  \subseteq \mathcal{N}(3)$. We are going to show that 
each element of $\mathcal{N}(3)$ belongs to $\mathcal{G}(3) $ by constructing 
explicit examples.  

Let $B_3$ be the dipole graph of genus $3$, where all the edges are assigned length one (see Example~\ref{eg:dipole}). By Example~\ref{eg:dipole:3}, there exist $p, q, r \in B_3$ with 
$\underline{n}_{K}(p) = (1, 2, 3)$, $\underline{n}_{K}(q) = (1, 2, 4)$ and $\underline{n}_{K}(r) = (1, 3, 5)$. 

It remains to see that $(1, 2, 5) \in \mathcal{G}(3) $. 

Let $K_4$ be the complete graph (this is the same as the wheel graph of genus $3$ in 
Example~\ref{eg:wheel}), where all the edges are assigned length one. Let $v \in K_4$ be a trivalent point. Since $K_4$ is non-hyperelliptic, we have $\rank(2(v)) = 0$. 
Since $K \sim 4(v)$, the tropical Riemann--Roch formula  gives $\rank(4(v)) = 2 + \rank(K - 4(v)) = 2$. 
Then $\rank(3(v)) = 1$. It follows that 
$\underline{n}_{K}(v) = (1, 2, 5)$. 

We conclude that 
$
\mathcal{G}(3)  = \mathcal{N}(3),  
$ 
and we obtain Theorem~\ref{thm:main:1}~(2). 

\subsection*{Genus \texorpdfstring{$4$}{4} case}
We define $\mathcal{S}(4)^\prime$ to be the subset of $\mathcal{S}(4)$ that excludes 
$(1, n_2, n_3, n_4)$ with $3 \leq n_2 <  n_3 < n_4 \leq 7$ and $(n_2, n_3, n_4) \neq (3, 5, 7)$. Then by 
Lemma~\ref{lemma:useful:2}, we have $\mathcal{G}(4) \subseteq \mathcal{S}(4)^\prime$. 
Explicitly, we have 
\begin{align*}
\mathcal{S}(4)^\prime  & = 
\left\{
(1, 2, 3, 4), (1, 2, 3, 5), (1, 2, 3, 6), (1, 2, 3, 7), (1, 2, 4, 5), (1, 2, 4, 6), 
\right. 
\\
& 
\qquad 
\left. (1, 2, 4, 7), (1, 2, 5, 6), (1, 2, 5, 7), (1, 2, 6, 7), (1, 3, 5, 7) \right\}. 
\end{align*}
Comparing with Theorem~\ref{thm:main:1}~(3), we have to show that 
\begin{equation}
\label{eqn:contradict:Clifford}
(1, 2, 5, 7), \; (1, 2, 6, 7)
\end{equation}
are not of the form $\underline{n}_{K}(p)$ for a point $p$ on a tropical curve 
$\Gamma$ of genus $4$, and then the other elements of 
$\mathcal{S}(4)^\prime$ are of the form $\underline{n}_{K}(p)$ for a point $p$ on a tropical curve $\Gamma$ of genus $4$. 

\subsubsection*{Exclusion of $(1, 2, 5, 7)$ and $(1, 2, 6, 7)$}
To derive a contradiction, suppose that $(1, 2, 5, 7)$ (resp. $(1, 2, 6, 7)$) were 
of the form $\underline{n}_{K}(p)$ for a point $p$. 
It follows that $\rank((p)) = 0$, $\rank(2(p)) = 0$, $\rank(3(p)) = 1$, $\rank(4(p)) = 2$, $\rank(5(p)) = 2$ (resp. $\rank(5(p)) = 3$), $\rank(6(p)) = 3$, and $\rank(7(p)) = 3$.  

In either case, 
since $\rank(4(p)) = 2$, the tropical Riemann--Roch formula gives 
$2 = \rank(4(p)) = 1 - 4 + 4 + \rank(K-4(p))$, so $\rank(K-4(p)) = 1$. 
Thus $4(p)$ is a special divisor. By Clifford's Theorem~\ref{thm:Clifford}, we would then get that $\Gamma$ is a hyperelliptic curve, and by Lemma~\ref{lem:cor:Clifford}, $4(p)$ would be linearly equivalent to $4 (v_0)$. Since 
$\rank(2(p)) = 0$, we have $2(p) \not\sim 2 (v_0)$, thus 
$6(p) = 2(p) + 4(p)\not\sim 2(v_0) + 4 (v_0) = 6 (v_0)$. Then 
with the notation in Proposition~\ref{prop:rank:hyp:ell:KY}, 
we have $p_{\Gamma}(6(p)) = 2$ and $\deg(6(p)) - p_{\Gamma}(6(p)) = 4= g$. 
Then Proposition~\ref{prop:rank:hyp:ell:KY} 
implies that $\rank(6(p)) = 2$, which contradicts with~$\rank(6(p)) = 3$. 
We conclude that $(1, 2, 5, 7)$ and $(1, 2, 6, 7)$ do not belong to $\mathcal{G}(4)$. 

\subsubsection*{Using the dipole graph $B_4$}
Let $B_4$ be the dipole graph of genus $4$ with two vertices $v, v^\prime$ joined by $5$ edges, where all the edges are assigned length one (see Example~\ref{eg:dipole}). Let $e_0$ be an edge of $B_4$, which we identify with 
the closed interval $[0, 1]$. Then we have the following. 

\begin{Proposition}
\label{prop:banana}
On the tropical curve $B_4$, 
Weierstrass gap sequences at $p \in [0, 1] = e_0$ 
are given as in Table~\ref{table:B4}. 

\begin{table}[htb!]
\begin{tabular}{l|lllllllllll}
\text{$\rank(k(p))$ for $1 \leq k \leq 7$} & $1$ & $2$ & $3$ & $4$ & $5$ & $6$ & $7$ &  & $\underline{n}_{K}(p)$ & & $\wt(p)$ \\ \hline
\smallskip
$p \in \left[0, \frac14\right) \cup \left(\frac34, 1\right]$  & $0$ & $0$ & $0$ 
& $0$ & $1$  & $2$ & $3$  & & $(1, 2, 3, 4)$ && $0$  \\ \hline
\smallskip
$p \in \left[\frac14, \frac13\right) \cup \left(\frac23, \frac34\right]
$ & $0$ & $0$ & $0$  
& $1$ & $1$  & $2$ & $3$  & & $(1, 2, 3, 5)$ && $1$ \\ \hline
\smallskip
$p \in  \left[\frac13, \frac25\right) \cup 
\left(\frac35, \frac23\right]
$ & $0$ & $0$ & $1$  
& $1$ & $1$  & $2$ & $3$  & & $(1, 2, 4, 5)$ && $2$ \\ \hline
\smallskip
$p \in \left[\frac25, \frac35\right] \setminus \left\{\frac12\right\}$ & $0$ & $0$ & $1$  
& $1$ & $2$  & $2$ & $3$  & & $(1, 2, 4, 6)$ && $3$ \\ \hline
\smallskip
$p = \frac12$ & $0$ & $1$ & $1$  
& $2$ & $2$  & $3$ & $3$  & & $(1, 3, 5, 7)$ && $6$ \\ \hline
\end{tabular}
\caption{Weierstrass gap sequences on $B_4$}
\label{table:B4}
\end{table}
\end{Proposition}

\Proof
By symmetry, we may assume that $p \in [0, 1/2]$. 
Let $v = 0, q = 1/2$ and $v^\prime = 1$. 

By the tropical Riemann-Roch formula, 
we have $\rank(7 (p)) = \deg(7 (p)) - g = 3$. 

\smallskip
{\bf Step 1.}\quad
If $p  \in [0, 1/4)$, then 
$4(p) \sim 3(0) + (4p)$. Since $3(0) + (4p)$ is $v$-reduced, we have 
$\rank(4(p)) = 0$. Then $\rank(k(p)) = 0$ for $1 \leq k \leq 3$. 
Since $\rank(7 (p)) = 3$, we have 
$\rank(5(p)) = 1$ and $\rank(6(p)) = 2$.  
Thus $\underline{n}_{K}(p) = (1, 2, 3, 4)$. 

\smallskip
{\bf Step 2.}\quad
Suppose that $p  \in [1/4, 1/3)$. Then 
$3(p) \sim 2(v) + (3p)$, and since $2(v) + (3p)$ is $v$-reduced, we have 
$\rank(3(p)) = 0$. On the other hand, we have 
$4(p) \sim 2 (q) + (v) + (4p-1)$ and 
$5(p) \sim 2 (q) + 2 (v) + (5p-1)$, so  Proposition~\ref{prop:rank:hyp:ell:KY} implies that $\rank(4(p)) = \rank(5(p)) = 1$. Since $\rank(7 (p)) = 3$, we have 
$\rank(6(p)) = 2$.  
Thus $\underline{n}_{K}(p) = (1, 2, 3, 5)$. 

\smallskip
{\bf Step 3.}\quad
Suppose that $p  \in [1/3, 2/5)$. We have $\rank((p)) = \rank(2(p)) = 0$. Since $3(p) \sim 2(q) + (3p-1)$, Proposition~\ref{prop:rank:hyp:ell:KY} implies that $\rank(3(p)) = 1$. Further, we have $5(p) \sim 3(q) + 
(v) + \left(5p - \frac{3}{2}\right)$, Proposition~\ref{prop:rank:hyp:ell:KY} implies that $\rank(5(p)) = 1$. Then $\rank(4(p)) = 1$. Since $\rank(7(p)) = 3$, 
we have $\rank(6(p)) = 2$. We conclude that 
$\underline{n}_{K}(p) = (1, 2, 4, 5)$. 

\smallskip
{\bf Step 4.}\quad
Suppose that $p  \in [2/5, 1/2)$. 
For $k \leq 5$, we have 
$k(p) \sim (k-1)(q) + \left(\frac12 - k(\frac12-p)\right)$ and 
$\frac12 - k(\frac12-p) \geq 0$. Proposition~\ref{prop:rank:hyp:ell:KY} implies that 
$\rank(k(p)) = \lfloor \frac{k-1}{2}\rfloor$ for $k \leq 5$. 
We have $6(p) \sim 4(q) + (v) + (6p-2)$ and 
$\frac25 \leq 6p-2 < 1$. Proposition~\ref{prop:rank:hyp:ell:KY} implies that $\rank(6(p)) = 2$. Thus we get 
$\underline{n}_{K}(p) = (1, 2, 4, 6)$. 

\smallskip
{\bf Step 5.}\quad
Suppose that $p = 1/2$. Then
Proposition~\ref{prop:rank:hyp:ell:KY} implies that 
$\rank(k(p)) = \lfloor \frac{k}{2}\rfloor$. 
It follows that $\underline{n}_{K}(p) = (1, 3, 5, 7)$. 
\QED

By Proposition~\ref{prop:banana} there exist $p, q, r, s, t \in B_4$ with 
$\underline{n}_{K}(p) = (1, 2, 3, 4)$, $\underline{n}_{K}(q) = (1, 2, 3, 5)$, $\underline{n}_{K}(r) = (1, 2, 4, 5)$, 
$\underline{n}_{K}(s) = (1, 2, 4, 6)$, and $\underline{n}_{K}(t) = (1, 3, 5, 7)$.  

\subsubsection*{Using the wheel graph}
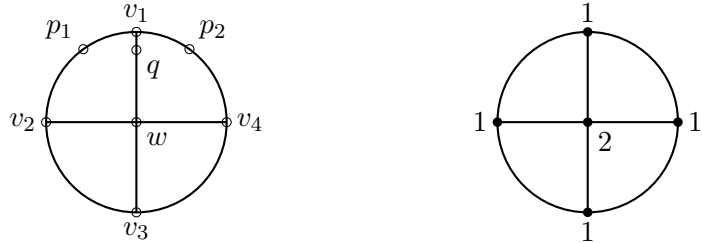
\begin{figure}[htb!]
\[
\begin{tikzpicture}
\begin{scope}[scale=0.6]
\draw [thick] (0,0) circle [radius=2];
\draw [thick] (0,0)--(0, 2);
\draw [thick] (0,0)--(-2, 0);
\draw [thick] (0,0)--(2, 0);
\draw [thick] (0,0)--(0, -2);
\draw (0, 0)circle (0.1);
\draw (0, 2)circle (0.1);
\draw (-2, 0)circle (0.1);
\draw (2, 0)circle (0.1);
\draw (0, -2)circle (0.1);
\draw (0, 0)node[below right]{$w$}; 
\draw (0, 2)node[above]{$v_1$}; 
\draw (-2, 0)node[left]{$v_2$}; 
\draw (2, 0)node[right]{$v_4$}; 
\draw (0, -2)node[below]{$v_3$}; 
\draw (0, 1.6)circle (0.1);
\draw (1.175570504, 1.618033989)circle (0.1);
\draw (-1.175570504, 1.618033989)circle (0.1);
\draw (0, 1.6)node[below right]{$q$}; 
\draw (1.175570504, 1.618033989)node[above right]{$p_2$}; 
\draw (-1.175570504, 1.618033989)node[above left]{$p_1$}; 
\end{scope}

\begin{scope}[shift={(6,0)}, scale=0.6]
\draw [thick] (0,0) circle [radius=2];
\draw [thick] (0,0)--(0, 2);
\draw [thick] (0,0)--(-2, 0);
\draw [thick] (0,0)--(2, 0);
\draw [thick] (0,0)--(0, -2);
\fill(0, 0) circle (3pt);
\fill(0, 2)circle (3pt);
\fill(-2, 0)circle (3pt);
\fill(2, 0)circle (3pt);
\fill(0, -2)circle (3pt);
\draw (0, 0)node[below right]{$2$}; 
\draw (0, 2)node[above]{$1$}; 
\draw (-2, 0)node[left]{$1$}; 
\draw (2, 0)node[right]{$1$}; 
\draw (0, -2)node[below]{$1$}; 
\end{scope}
\end{tikzpicture}
\]
\caption{The wheel graph of genus $4$ and its canonical divisor.}
\label{figure:wheel:genus:4}
\end{figure}

Let $\Gamma$ be the wheel graph of genus $4$, where all the edges are assigned length one (see Example~\ref{eg:wheel} and Figure~\ref{figure:wheel:genus:4}). 
We write $w$ for the central vertex (the center of the wheel), and 
$v_1, \ldots, v_4$ for the other vertices. 
We identify the edge $v_1w$ with the closed interval $[0, 1]$ with $v_1 = 0$ and $w = 1$. 
Let $q \in v_1w$ be the point with coordinate $q = 1/5$. 

\begin{Lemma}
\label{lem:wheel:4}
\begin{enumerate}
\item
$\underline{n}_{K}(w) = (1, 2, 3, 7)$. 
\item
$\underline{n}_{K}(q) = (1, 2, 3, 6)$. 
\end{enumerate}
\end{Lemma}

\Proof
We have 
$
K = 2 (w) + (v_1) + (v_2) + (v_3) + (v_4),
$
and we see that $K \sim 6 (w) \sim 4 (v_1) +  (w) + (v_3)$. 

(1) Since $3(w)$ is $v_1$-reduced, we have $\rank(3(w)) = 0$. 
By the tropical Riemann--Roch formula, we have $\rank(6(w)) = 1 - 4 + 6 + \rank(K - 6(w)) = 3$. Then we have $\rank(4(w)) = 1$ and  
$\rank(5(w)) = 2$. 
We obtain $\underline{n}_{K}(w) = (1, 2, 3, 7)$. 

(2) Identifying the edge $v_1v_2$ (resp. $v_1 v_4$) with $[0, 1]$, let 
$p_1$ (resp. $p_2$) be the point with coordinate $2/5$. Then we see that 
$3(q) \sim (p_1) + (p_2) + (w)$. Since $(p_1) + (p_2) + (w)$ is $v_3$-reduced, 
we have $\rank(3(q)) = 0$. Since $5(q) \sim 4 (v_1)+ (w)$, 
by the tropical Riemann--Roch formula, we have $\rank(5(q)) = 1 - 4 + 5 + \rank(K - 5(q)) = 
2+ \rank((v_3)) = 2$. It follows that $\rank(4(q))) = 1$. 
Further, $\rank(6(q)) = 1 - 4 + 6 + \rank(K - 6(q)) = 
3+ \rank(-(q) + (v_3)) = 2$. We obtain $\underline{n}_{K}(w) 
= (1, 2, 3, 6)$. 
\QED

By Lemma~\ref{lem:wheel:4}, there exist a tropical curve of $\Gamma$ of genus $4$ and $p, q \in \Gamma$ such that  $\underline{n}_{K}(p) = (1, 2, 3, 6)$ and $\underline{n}_{K}(q) = (1, 2, 3, 7)$.

It remains to show that 
\[
(1, 2, 4, 7), (1, 2, 5, 6)
\]
appear as a Weierstrass gap sequence of a point on a tropical curve. 

\subsubsection*{Construction for $(1, 2, 4, 7)$}
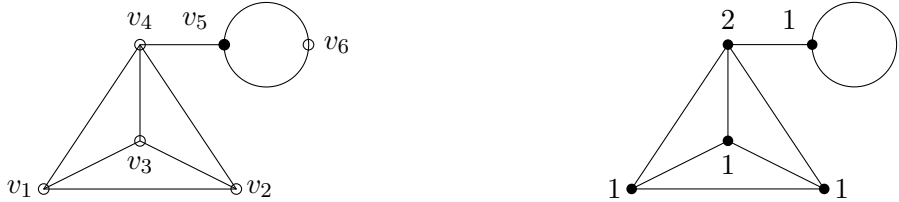
\begin{figure}[ht!]
\begin{minipage}{0.45\textwidth}
\centering
\begin{tikzpicture}[scale=.8]
	\coordinate (A) at (0,0);
	\coordinate (B) at (3.2,0);
	\coordinate (C) at (1.6,2.4);
	\coordinate (D) at (1.6,0.8);
	\coordinate (E) at (3, 2.4);
	\coordinate (F) at (4.4, 2.4);
	
	\draw (A) -- (B) -- (C) -- cycle;
	\draw (A) -- (D) -- (B);
	\draw (C) -- (D);
	\draw (C) -- (E);
	
	\node[left] at (A) {$v_1$};
	\node[right] at (B) {$v_2$};
	\node[below=2pt] at (D) {$v_3$};
	\node[above=2pt] at (C) {$v_4$};
	\node[right=2pt] at (F) {$v_6$};

	\draw (3.7, 2.4) circle (20pt);
	
	\fill  (E) circle (2.5pt); \node[above left=2pt] at (E) {$ v_5$}; 	
	\foreach \c in {A,B,C,D, E, F} {
\draw (\c) circle (2.5pt);
	}
	
\end{tikzpicture}
\end{minipage}
\begin{minipage}{0.45\textwidth}
\centering
\begin{tikzpicture}[scale=.8]
	\coordinate (A) at (0,0);
	\coordinate (B) at (3.2,0);
	\coordinate (C) at (1.6,2.4);
	\coordinate (D) at (1.6,0.8);
	\coordinate (E) at (3, 2.4);
	
	\draw (A) -- (B) -- (C) -- cycle;
	\draw (A) -- (D) -- (B);
	\draw (C) -- (D);
	\draw (C) -- (E);
	
	\draw (3.7, 2.4) circle (20pt);
	
	\node[left] at (A) {$1$};
	\node[right] at (B) {$1$};
	\node[below=2pt] at (D) {$1$};
	\node[above=2pt] at (C) {$2$};
	\node[above left=2pt] at (E) {$1$};
	
	\foreach \c in {A,B,C,D, E} {
\fill (\c) circle (2.5pt); 
	}
	
\end{tikzpicture}
\end{minipage}
\caption{The tropical curve $\Gamma$ and its canonical divisor.}
\label{fig:complete4circle}
\end{figure}

We consider the tropical curve $\Gamma$ of genus $4$ with 
the set of vertices $\{v_1, \ldots, v_6\}$ as in Figure~\ref{fig:complete4circle}, where all the edges are assigned length one. 

Let $p$ be the vertex $v_5$. Since $\Gamma$ is neither a tree nor a hyperelliptic graph, we have 
$\rank((p)) = \rank(2 (p)) = 0$. 
Since $3(p) \sim 3 (v_4) \sim (v_1) + (v_2) + (v_3)$ and $3(p) \sim (p) + 2 [v_6]$, 
we have $\rank(3(p)) = 1$. Note that 
$K = (v_1) + (v_2) + (v_3) + 2 (v_4) +  (v_5) \sim 6 (p)$. For $k = 4, 5, 6$, 
the tropical Riemann-Roch formula implies $\rank(k (p)) = (k-3) + \rank((6-k)(p))$. Then $\rank(4(p)) = 1$, $\rank(5(p)) = 2$, and $\rank_{\Gamma_1}(6(p)) = 3$. It follows that $\underline{n}_{K}(p) = 
(1, 2, 4, 7)$. 

\subsubsection*{Construction for $(1, 2, 5, 6)$} 
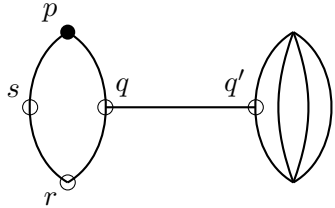
\begin{figure}[htb]
\[
\begin{tikzpicture}
\draw[thick] (0,0) to [out=60,in=120,relative] (0, 2);
\draw[thick] (0,0) to [out=-60,in=-120,relative] (0, 2);
\draw[thick] (0.5,1)--(2.5,1); 
\draw[thick] (3,0) to [out=60,in=120,relative] (3, 2);
\draw[thick] (3,0) to [out=-60,in=-120,relative] (3, 2);
\draw[thick] (3,0) to [out=20,in=160,relative] (3, 2);
\draw[thick] (3,0) to [out=-20,in=-160,relative] (3, 2);
\draw (0,2)node[above left]{$p$}; 
\fill[black] (0, 2)circle (0.1);
\draw (0.5, 1)node[above right]{$q$}; 
\draw (2.5, 1)node[above left]{$q^\prime$}; 
\draw (0, 0)node[below left]{$r$}; 
\draw (0, 0)circle (0.1); 
\draw (2.5, 1)circle (0.1); 
\draw (0.5, 1)circle (0.1); 
\draw (-0.5, 1)circle (0.1); 
\draw (-0.5, 1)node[above left]{$s$}; 
\end{tikzpicture}
\]
\caption{The tropical curve $B_1$-$B_3$ joining $B_1$ and $B_3$}
\label{fig:B1:B3}
\end{figure}

Let $B_3$ be the dipole graph of genus $3$ (edge lengths can be arbitrary), and let $B_1$ be the circle (of an arbitrary length); see Example~\ref{eg:dipole}. 
Take $p, r$ as antipodal points on $B_1$, and let 
$q$ be the middle point of an edge connecting $p$ and $r$ on $B_1$. 
Let $q'$ be the middle point of an edge of $B_3$. 
We write for $B_1$-$B_3$ for the tropical hyperelliptic curve joining $B_1$ and $B_3$ by a bridge (of an arbitrary length) with end points $q, q'$. 
Let $s$ be the antipodal point of $q$ in $B_1$. 
We consider the tropical curve $\Gamma = \text{$B_1$-$B_3$}$, which has genus $4$; see Figure~\ref{fig:B1:B3}. 

Let us compute $\rank(k(p))$ for $1 \leq k \leq 6$. 
We have $\rank((p)) = 0, \rank(2(p)) = 0$. For $k = 3, 4, 5, 6$, we have 
\[
3(p) \sim  2(q) + (r), \quad 4(p) \sim 4(q), \quad  5(p) \sim  4(q) + (p), \quad 
6(p) \sim  5(q) + (s). 
\]
since $2(q) + (r), 4(q), 4(q) + (p), 5(q) + (s)$ are all $q$-reduced, Proposition~\ref{prop:rank:hyp:ell:KY} gives 
\[
\rank(3(p)) = 1, \quad 
\rank(4(p)) = 2, \quad 
\rank(5(p)) = 2, \quad 
\rank(6(p)) = 2. 
\]
We obtain $\underline{n}_{K}(p) = (1, 2, 5, 6)$. 

\smallskip
This completes the proof of Theorem~\ref{thm:main:1}~(3). 
\QED

\begin{Remark}
\label{rmk:no:numerical:semigp} 
We have used $B_4$ to construct a point $p$ with $\underline{n}(p) = (1, 2, 4, 6)$, 
for which the complement in $\ZZ_{>0}$ is not a numerical subgroup.  

This can be generalized for any genus $g \geq 4$. Indeed, let $B_g$ be the dipole graph of genus $g$, where all the edges are assigned length one (see Example~\ref{eg:dipole}). With the notation there, we identify the edge $e_0$ with 
the closed interval $[0, 1]$ with $v = 0$ and $v^\prime = 1$. Let $q = 1/2$ be the middle point of $e_0$. 
We take  $0 < \varepsilon < \frac{1}{2(2g-3)}$ and put $p_1 = 1/2 - \varepsilon \in e_0$. 

By Proposition~\ref{prop:rank:hyp:ell:KY}, we have 
$
p_{B_g}(k (p_1)) = \lfloor \frac{k-1}{2} \rfloor 
$
for $1 \leq k \leq 2g-3$. Since $\deg(k (p_1)) - \lfloor \frac{k-1}{2} \rfloor = \lceil \frac{k+1}{2} \rceil  \leq 
g - 1 \; (< g)$ for $1 \leq k \leq 2g-3$, we have  
$
\rank(k (p_1)) 
= 
p_{B_g}(k (p_1))$. 
Further, since $(2g-2)(p_1) \sim 
(2g-4)(q) + (v) + (2p_1)$, Proposition~\ref{prop:rank:hyp:ell:KY}
implies that $\rank((2g-2) (p_1)) = g-2$. 
Thus 
$
\underline{n}_{K}(p_1) 
= (1, 2, 4, 6, \ldots, 2g-2) 
$. 
Since $g \geq 4$, we have $6 \in \underline{n}_{K}(p_1)$, and $3 \not\in \underline{n}_{K}(p_1)$. 
Thus ${\ZZ_{> 0}} \setminus \underline{n}_{K}(p_1)$ is not a numerical semigroup. 
\end{Remark}

\section{Questions}
\label{sec:open:q}
We end the paper by posing several questions. 
Let $\Gamma$ be a tropical curve of genus $g\geq 2$. 
Let $D$ be a divisor on $\Gamma$ of degree $d$ and nonnegative rank $r$, and let $B_1, \ldots, B_m$ be the set of all maximal $D$-Weierstrass loci of $\Gamma$ (for the definition, see the discussion after Proposition~\ref{prop:specialization:coordinatewise} in Section~\ref{sec:intro}). We set $\wt_D(B_i) \colonequals \wt_D(p_i)$ for any $p_i \in B_i$, which does not depend on the choice of $p_i$. 

\begin{Question}
\label{q:max:W:1}
Suppose that $D$ is the canonical divisor $K$. 
Can one give a lower and an upper bound of the total sum 
\[
\sum_{i=1}^m \wt_K(B_i)? 
\]
\end{Question}
The computations in Section~\ref{sec:proof:thm:max:WP:genus:3} show that, when $g = 2$, 
one has $\sum_{i=1}^m \wt_K(B_i) = 6$ and that when $g = 3$, 
one has $8 \leq \sum_{i=1}^m \wt_K(B_i) \leq 24$, which are 
the best possible lower and upper bounds.  By Proposition~\ref{prop:contrast} and Remark~\ref{rmk:uniformity}, we know that there is a constant $c(g)$ depending only on $g$ such that $\sum_{i=1}^m \wt_K(B_i) \leq c(g)$.

\begin{Question}
\label{q:B:wt}
Let $A$ be a $D$-Weierstrass locus, and let $B_1, \ldots, B_k$ be the set of maximal $D$-Weierstrass loci that are included in $A$. Can one compare $\sum_{i=1}^k \mu_D(B_i)$ and $\mu_D(A)$? Can one compare $\mu_D(B_i)$ and $\wt_D(B_i)$?  
\end{Question}

We show in Corollary~\ref{cor:M:m} that for any $b$ (not necessarily contained in a maximal $D$-Weierstrass locus), we have $\mu_D(\{b\}) \leq \wt_D(b)$. 

\begin{Question}
\label{q:bridgeless}
Assume that $\Gamma$ is bridgeless. Can one give an example of $\Gamma$ such that a maximal Weierstrass locus $B$ (for the canonical divisor $K$) is not a singleton?
\end{Question}
For $n = 2, 3$, the computations in Section~\ref{sec:proof:thm:max:WP:genus:3} show 
that any maximal Weierstrass locus is a singleton. 

\begin{Question}
For $\underline{n} \in \mathcal{S}(r, d)$, 
let $B$ be a connected component of $\WL_{\geq  \underline{n}}(D)$. Can one define a refined weight 
$\widetilde{\mu}_D(B)$ and improve the inequality 
in Theorem~\ref{thm:M:m}? 
\end{Question}

\begin{Question}
\label{q:max:W:3}
Put 
\begin{align*}
\mathcal{MG}(g)
& \colonequals  \left\{\underline{n}\in\mathcal{S}(g)  \;\left|\; 
\begin{aligned}
& \text{$\underline{n} = \underline{n}_{K}(p)$ for some tropical curve $\Gamma$ of genus $g$} 
\\
& \text{and a {\em maximal} Weierstrass point $p \in \Gamma$}
\end{aligned}
\right.  \right\}, 
\\
\widetilde{\mathcal{G}}(g)
& \colonequals  \left\{\underline{n}\in\mathcal{S}(g)  \;\left|\; 
\begin{aligned}
& \text{$\underline{n} = \underline{n}_{K_{\mathcal{C}}}(x)$ for some smooth projective curve $\mathcal{C}$ of genus $g$ over} 
\\
& \text{an algebraically closed field $k$ and a Weierstrass point $x \in \mathcal{C}(k)$}
\end{aligned}
\right.  \right\}. 
\end{align*}
Then $\mathcal{MG}(g) \subseteq \mathcal{G}(g)$. 
Can one compare $\widetilde{\mathcal{G}}(g)$, $\mathcal{MG}(g)$, and $\mathcal{G}(g)$ in a way that contributes to studying  $\mathcal{MG}(g)$ or $\mathcal{G}(g)$? 
\end{Question}

\bibliographystyle{alpha}
\bibliography{bibliography}
\end{document}